\documentclass[a4paper,oneside,12pt]{article}
\usepackage{amsmath,amsfonts,amssymb,amsthm,mathrsfs}
\usepackage[a4paper,vmargin={3.5cm,3.5cm},hmargin={2.5cm,2.5cm}]{geometry}
\usepackage[font=sf, labelfont={sf,bf}, margin=1cm]{caption}
\usepackage{graphicx,graphics}
\usepackage{epsfig}
\usepackage{latexsym}
\usepackage[applemac]{inputenc}
\usepackage{ae,aecompl}
\usepackage[english]{babel}
 \usepackage[colorlinks=true]{hyperref}
\usepackage{pstricks}

\newtheorem{theorem}{Theorem}

\newtheorem{open}{Open Question}
\newtheorem{proposition}{Proposition}

\newtheorem{lemma}{Lemma}
\theoremstyle{definition}

\theoremstyle{remark}
\newtheorem{remark}{Remark}

\newcommand{\R}{\mathbb R}

\newcommand{\Z}{\mathbb Z}

\renewcommand{\geq}{\geqslant}
\renewcommand{\leq}{\leqslant}

\renewcommand{\leq}{\leqslant}
\renewcommand{\geq}{\geqslant}

 \title{\bf Uniform infinite planar quadrangulations with a boundary}
  \author{Nicolas Curien and Grégory Miermont}
  \date{
  }
\begin{document}
\maketitle

\abstract{We introduce and study the uniform infinite planar quadrangulation (UIPQ) with a boundary via an extension of the construction of \cite{CMMinfini}. We then relate this object to its   simple boundary analog using a pruning procedure. This enables us to study the aperture of these maps, that is, the maximal distance between two points on the boundary, which in turn sheds new light on the geometry of the UIPQ. In particular we prove that the self-avoiding walk on the UIPQ is diffusive.}


\section*{Introduction} Motivated by the theory of $2$D quantum gravity, the probabilistic theory of random planar maps has been considerably growing over the last few years. In this paper we continue the study of the geometry of random maps and focus in particular on random quadrangulations with a \emph{boundary}. \medskip 

Recall that a planar map is a proper  embedding 
of a finite connected planar graph into the two-dimensional sphere seen up to 
orientation-preserving homeomorphisms. The faces are the connected components of the complement of the union of the edges, and the degree of a face is the number of edges that are incident to it, where it should be understood that an edge is counted twice if it lies entirely in the face.  A map is a quadrangulation if all its faces have degree $4$.  All the maps considered in this work are rooted, meaning that an oriented edge is distinguished and called the {\em root edge}. The face lying to the right of the root edge is called the {\em root face}. 
 \begin{figure}[!h]
 \begin{center}
 \includegraphics[width=13cm]{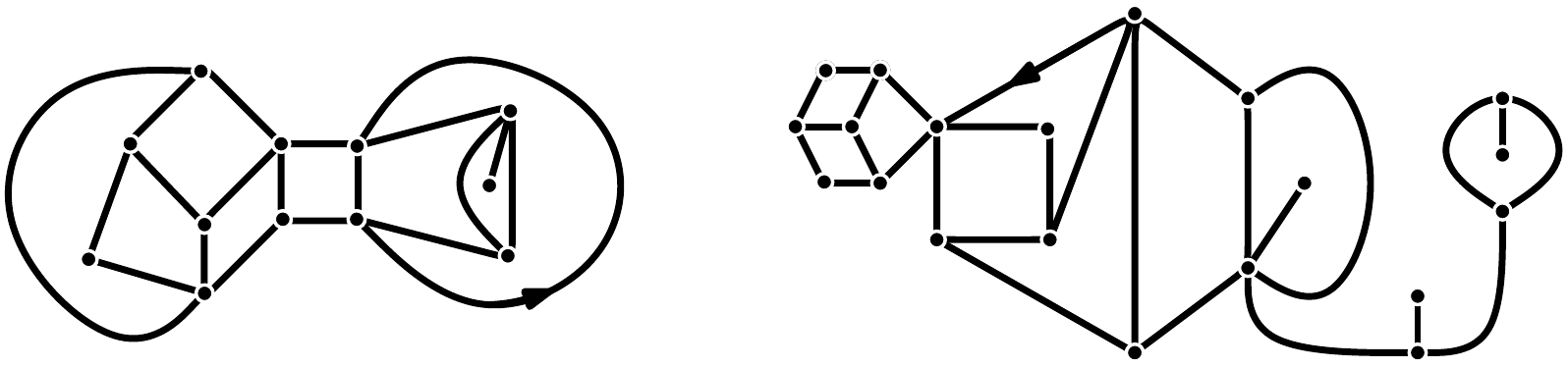}
 \caption{A quadrangulation with a simple and a non-simple boundary.}
 \end{center}
 \end{figure}
 
A planar map $q$ is a \emph{quadrangulation with a boundary} if all its faces have degree four except possibly the root face, which can have an arbitrary even degree. Since we want to consider this distinguished face as lying ``outside'' of the map, we also call it the \emph{external face}. The degree (which must be even) of the external face is called the {\em perimeter} of the map, and the boundary is said to be simple if during the contour of the external face all the vertices on the boundary are visited only once (i.e. there are no pinch points on the boundary). The size of $q$ is its number of faces minus one.


 Uniform quadrangulations of size $n$ with a boundary of perimeter $p$ have recently been studied  from a combinatorial and a probabilistic point of view \cite{Bet11,BG09}. Three different regimes have to be distinguished: If $p \ll n^{1/2}$, then these maps converge, in the scaling limit, towards the Brownian map introduced in \cite{LG11,Mie11}. If $ p \asymp n^{1/2}$, then the boundary becomes macroscopic and Bettinelli \cite{Bet11} introduced the natural candidate for the scaling limits of these objects which is a sort of Brownian map with a hole. When $ p \gg n^{1/2}$ these random quadrangulations fold on themselves and become tree-like \cite{Bet11,BG09}. In this work, we shall take a different approach and study infinite local limits of quadrangulations with a boundary as the size tends to infinity. Let us precise the setting. \medskip 
 
 In a pioneering work \cite{BS01}, Benjamini \& Schramm initiated the study of local limits of maps.  If $m,m'$ are two rooted maps, the local distance between $m$ and $m'$ is
 \begin{eqnarray*} \mathrm{d_{map}}(m,m') &=& \big(1+ \sup\{r\geq 0 : B_{r}(m)=B_{r}(m')\}\big)^{-1}, \end{eqnarray*}
where $B_{r}(m)$ denotes the map formed by the faces of $m$ whose vertices are \emph{all} at graph distance smaller than or equal to $r$ from the origin of the root edge in $m$. Let $Q_n$ be  uniformly distributed over the set of all rooted quadrangulations with $n$ faces. Krikun \cite{Kri05} proved that
\begin{eqnarray} Q_n &\xrightarrow[n\to \infty]{(d)}& Q_\infty, \label{def:UIPQ}\end{eqnarray} in distribution in the sense of $ \mathrm{d_{map}}$. The object $Q_\infty$  is a random infinite rooted quadrangulation called the Uniform Infinite Planar Quadrangulations (UIPQ) (see also \cite{Ang03,AS03} for previous works concerning triangulations). The UIPQ (and its related triangulation analog) has been the subject of numerous researches  in recent years, see \cite{CD06,CMMinfini,Kri05,Kri08,LGM10}. Despite these progresses, the geometry of the UIPQ remains quite mysterious. The purpose of this article is to provide some new geometric understanding of the UIPQ via the study of UIPQ \emph{with a boundary}.\\
 
We will show that the convergence \eqref{def:UIPQ} can be extended to quadrangulations with a boundary.  More precisely, for any $p \geq 1$, we let $Q_{n,p}$ (resp.\, $ \widetilde{Q}_{n,p}$) be a uniform quadrangulation of size $n$ and with a (resp.\,simple) boundary of perimeter $2p$ then we have 
 \begin{eqnarray*} 
  \widetilde{Q}_{n,p} & \xrightarrow[n\to\infty]{(d)} & \widetilde{Q}_{\infty,p},\\
  Q_{n,p} & \xrightarrow[n\to\infty]{(d)} & Q_{\infty,p},  \end{eqnarray*}
in distribution for the metric $ \mathrm{d_{map}}$. The random maps $Q_{\infty,p}$ and $\widetilde{Q}_{\infty,p}$ are called the Uniform Infinite Planar Quadrangulation with a (resp.\,simple) boundary of perimeter $2p$. The first convergence is an easy consequence of \eqref{def:UIPQ} (see the discussion around \eqref{def:UIPQs} below) whereas the second convergence requires an adaptation of the techniques of \cite{CMMinfini}: In Theorem \ref{CMM10bound}, we construct $Q_{\infty,p}$ from a labeled ``treed bridge" and extend the main result of \cite{CMMinfini} to our setting. This construction is yet another example of the power of the bijective technique triggered by Schaeffer \cite{Sch98} which has been one of the key tool for studying random planar maps \cite{BDFG04,CMS09,CD06,Mie09}.

\begin{figure}[!h]
 \begin{center}
 \includegraphics[width=9cm]{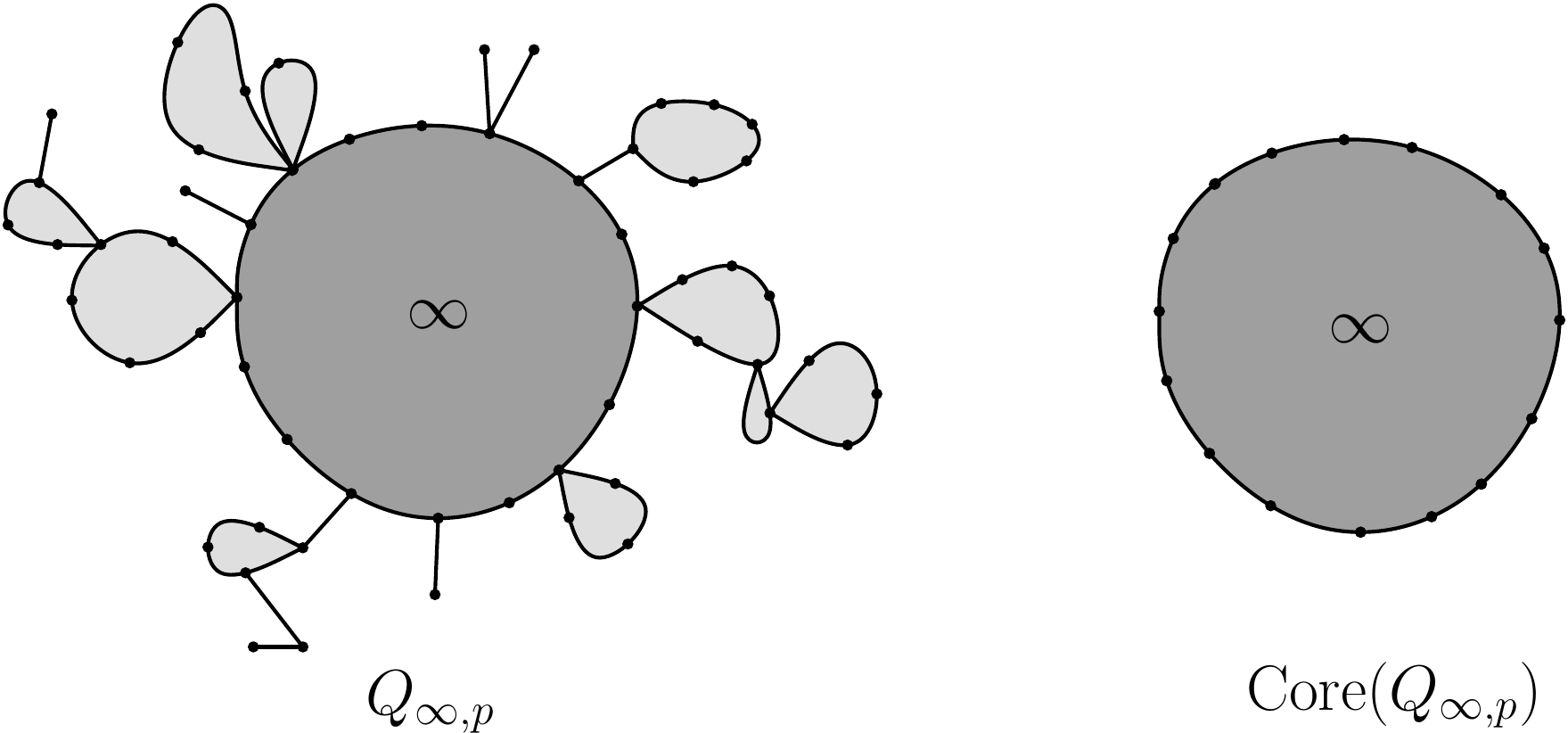}
 \caption{ \label{pruning-intro} Pruning of $Q_{\infty,p}$.}
 \end{center}
 \end{figure}

We then turn to the study of the UIPQ's with a boundary and their relationships. Although well-suited for the definition and the study of $Q_{\infty,p}$, the techniques ``à la Schaeffer" seem much harder to develop in the case of simple boundary because of the topological constraint imposed on the external face.  In order to bypass this difficulty we use a pruning decomposition to go from non-simple to simple boundaries. More precisely, we  prove that $Q_{\infty,p}$ has a unique infinite irreducible component, that is, a core made of an infinite quadrangulation with a simple boundary together with finite quadrangulations hanging off from this core, see Fig.\,\ref{pruning-intro}.
 
  We show that if we remove these finite components then the core of $Q_{\infty,p}$ has a  (random) perimeter $|\partial \mathrm{Core}(Q_{\infty,p})|$ which is roughly a third of the original one. More precisely we prove in Proposition \ref{1/3} the following convergence in distribution 
   \begin{eqnarray*}\frac{|\partial \mathrm{Core}(Q_{\infty,p})|-2p/3}{p^{2/3}} & \xrightarrow[p\to\infty]{(d)}& \mathcal{Z},  \end{eqnarray*} where $ \mathcal{Z}$ is a spectrally negative stable random variable of parameter $3/2$. Furthermore, conditionally on its perimeter, the core is  distributed as a UIPQ with a simple boundary (Theorem \ref{thm:fixed}). This confirms and sharpens a phenomenon already observed in a slightly different context by Bouttier \& Guitter, see \cite[Section 5]{BG09}

\bigskip As an application of our techniques we study the \emph{aperture} of these maps: If $q$ is a quadrangulation with a boundary, we denote the maximal graph distance between two vertices on the boundary of $q$ by $\operatorname{aper}(q)$ and call it the {aperture} of $q$. We prove that the aperture of the UIPQ with a simple boundary of perimeter $p$ is strongly concentrated around $\sqrt{p}$, in the sense of the following statement. 
 
  \begin{theorem} \label{main}
There exists $c,c'>0$ such that for all $p\geq 1$ and $\lambda >0$ the  aperture of a uniform infinite planar  quadrangulation with simple boundary of perimeter $2p$ satisfies
\begin{eqnarray*} P\left( \mathrm{aper}(\widetilde{Q}_{\infty,p}) \geq  \lambda \sqrt{p} \right)  &\leq& c \,  p^{2/3}\, \exp\big(-c'\lambda^{2/3}\big) . 
\end{eqnarray*}
\end{theorem}
This result is first established for UIPQ with general boundary using the construction from a treed bridge (Proposition \ref{diamext}) and then transferred to the simple boundary case using the pruning procedure.  This theorem provides a new tool for studying the UIPQ itself via the technique of \emph{peeling}, see \cite{Ang03,BCsubdiffusive}. In particular, Theorem \ref{main} is one of the key estimates of \cite{BCsubdiffusive} used to prove that the simple random walk on the UIPQ is subdiffusive with exponent less than $1/3$.\medskip

\begin{figure}[!h]
 \begin{center}
 \includegraphics[height=5cm]{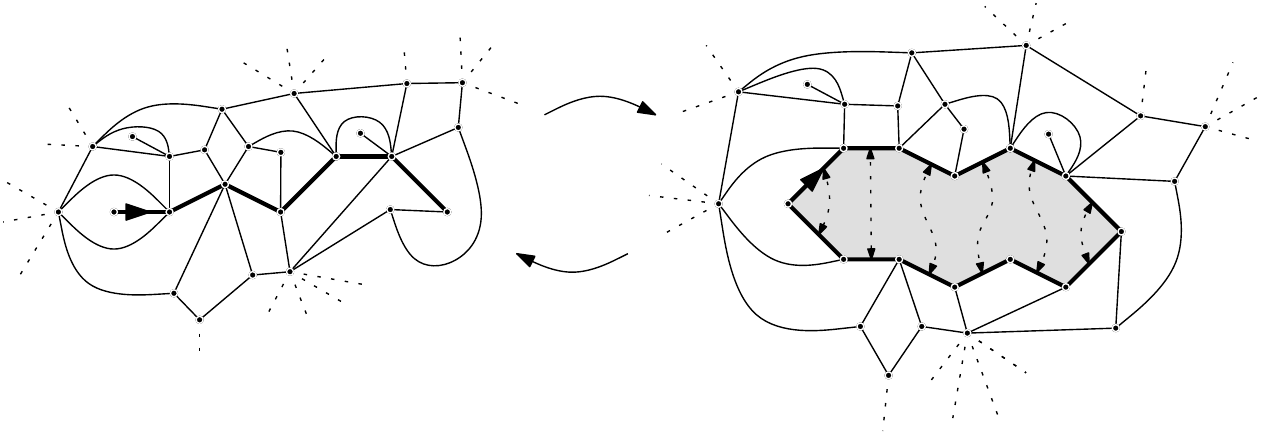}
 \caption{ \label{zipper}Zip the external face and unzip the self-avoiding walk.}
 \end{center}
 \end{figure}

Let us finish this introduction with one more motivation. 
There is an obvious bijective correspondence between, on the one hand, quadrangulations of size $n$ with a self-avoiding path of length $p$ starting at the root edge and, on the other hand, quadrangulations with simple boundary of perimeter $2p$ and size $n$: Simply consider the self-avoiding walk as a zipper. See Fig.\,\ref{zipper}.  Hence, the UIPQ with simple boundary of perimeter $2p$ can be seen as an annealed model of UIPQ endowed with a self-avoiding path of length $p$.

 With this correspondence, the aperture of the map obviously bounds the maximal graph distance of any point of the self-avoiding walk to the origin of the map. The estimates of Theorem \ref{main} then show that, when $p$ is large, the maximal graph distance displacement of the self-avoiding walk with respect to the root of the map is at most of order $\sqrt{p}$. This contrasts with the Euclidean case where a displacement of order $p^{3/4}$ is conjectured.

 Let us remark, however, that the aperture of a quadrangulation with boundary only gives an upper bound on the maximal displacement of the SAW obtained after zipping.
 \begin{open} Consider the infinite quadrangulation with a self-avoiding walk obtained after zipping the boundary of $ \widetilde{Q}_{\infty,p}$.  Prove a lower bound (if possible matching the order $\sqrt{p}$) on the maximal displacement from the root of this walk as $p \to \infty$.
 \end{open}

The paper is organized as follows: The first section contains some background on quadrangulations with a boundary. In the second section, we present the bijective techniques adapted from \cite{BG09} that we apply in Section \ref{sec:UIPQbound} to define the UIPQ with general boundary and study its  aperture. The fourth section is devoted to the pruning decomposition and its applications. Finally, the last section contains applications, extensions and comments. In particular, we define the UIPQ of the half-plane (with infinite boundary) with general and simple boundary and propose some open questions.\\

\noindent \textbf{Acknowledgments: } We are grateful to Jérémie Bettinelli for a careful reading and numerous comments on a first version of this article.

\section{Quadrangulations with a boundary}
\subsection{Definitions}

Recall that all the maps we consider are \emph{rooted}, that is given with one distinguished oriented edge $\vec{e}$.

 A planar map $q$ is a \emph{quadrangulation with a boundary} if all its faces have degree four, with the possible exception of the root face (also called external face). Since quadrangulations are bipartite, the degree of the external face has to be even. The boundary of the external face is denoted by $\partial q$ and its degree by $|\partial q|$. We say that $q$ has a perimeter $|\partial q|$ and its \emph{size} $|q|$ is the number of faces minus $1$.

A quadrangulation has a \emph{simple} boundary if there is no pinch point on the boundary, that is, if $\partial q$ is a cycle with no self intersection. By convention, all the notation involving a simple quadrangulation will be decorated with a ``$\sim$'' to avoid confusion.

We denote by $\mathcal{Q}_{n,p}$ (resp.\,$\widetilde{ \mathcal{Q}}_{n,p}$) the set of all rooted quadrangulations with (resp.\,simple) boundary with $n+1$ faces and such that the external face has degree $2p$ and by $q_{n,p}$ (resp.\,$\widetilde{q}_{n,p}$) its cardinal. By convention, the set $\mathcal{Q}_{0,0}= \widetilde{ \mathcal{Q}}_{0,0}$ contains a unique ``vertex'' map denoted by $\dag$. Note also that $\mathcal{Q}_{0,1}= \widetilde{ \mathcal{Q}}_{0,1}$ is composed of the map with one oriented edge (which has simple boundary). Note that any quadrangulation with boundary of perimeter $2$ can be seen as a rooted quadrangulation without boundary, by contracting the external face of degree two.

\subsection{Enumeration} \label{warmup} 
\begin{figure}[!h]
 \begin{center}
 \includegraphics[height=4cm]{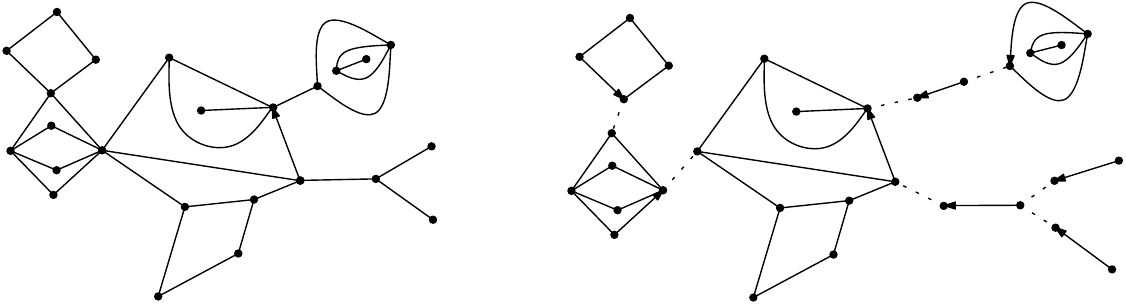}
 \caption{ \label{fig:decomposition} Decomposition of a quadrangulation with boundary into irreducible components.}
 \end{center}
 \end{figure}

Let $q$ be a quadrangulation with boundary. If the boundary of $q$ is not simple (if it has some separating vertices) we can decompose $q$ unambiguously into quadrangulations with simple boundary attached by the separating vertices of the boundary of $q$: These quadrangulations are called \emph{the irreducible components} of $q$. The root edge is carried by a unique irreducible component, and all other irreducible components have a unique boundary vertex which is closest to the component of the root. By convention, we root each component at the oriented edge that immediately precedes this particular vertex in counterclockwise order. See Fig.\ref{fig:decomposition}.

%

 We gather here a few enumeration results that will be useful in the following. We refer to \cite{BG09} for the derivations of these formul\ae.\medskip

If $q$ is a quadrangulation with a general boundary we can also decompose $q$ according to the irreducible component that contains its root edge and other quadrangulations with general boundary attached to it. This decomposition yields an identity relating the bi-variate generating functions of quadrangulations with simple and general boundary: For $g,z \geq 0$,  let $W$ (resp.\,$\widetilde{W}$) be the bi-variate generating function of $\mathcal{Q}_{n,p}$ (resp.\,$\widetilde{\mathcal{Q}}_{n,p}$) with weight $g$ per internal face and $\sqrt{z}$ per edge on the boundary,  that is
 \begin{eqnarray*} W(g,z) := \sum_{n,p \geq 0}{q}_{n,p}g^n z^p, \quad \mbox{ and } \quad \widetilde{W}(g,z) := \sum_{n,p \geq 0} {\widetilde{q}}_{n,p}g^n z^p.  \end{eqnarray*}
Then the last decomposition  translates into the identity   \begin{eqnarray} \label{identityGF}
\widetilde{W}\big(g, z W^2(g,z)\big) &=& W(g,z).  \end{eqnarray}
The exact expression of $W$ can be found in \cite{BG09}, it reads 

 \begin{eqnarray} W(g,z) =  \omega(1-gR^2(\omega-1)),  \ 
\mbox{ where } \ \omega = \frac{1-\sqrt{1-4z R}}{2zR} \ \mbox{ and } \  R = \frac{1-\sqrt{1-12g}}{6g}. \label{WWW}\end{eqnarray}
From this we can deduce 
\begin{eqnarray*}
{q}_{n,p} &=& \frac{(2p)!}{p!(p-1)!}3^n\frac{(2n+p-1)!}{n!(n+p+1)!},\end{eqnarray*} for $n \geq 0$ and $p\geq 1$.
Note the asymptotics 
\begin{eqnarray}
 q_{n,p} & \underset{n \to \infty}{\sim}& C_{p} 12^n n^{-5/2}, \label{asympqnp}\\
C_{p} &=& \frac{(2p)!}{p!(p-1)!}2^{p-1} \pi^{-1/2} \underset{p \to \infty}{\sim} (2\pi)^{-1} 8^p\sqrt{p} \nonumber.
\end{eqnarray}Moreover equation \eqref{identityGF} enables us to find the expressions for $\widetilde{q}_{n,p}$ (see \cite{BG09} for more details) namely
\begin{eqnarray}
\widetilde{q}_{n,p} & = & 3^{-p} \frac{(3p)!}{p!(2p-1)!} 3^n\frac{(2n+p-1)!}{(n-p+1)!(n+2p)!}, \quad \mbox{ for }n \geq 1 \mbox{ and } p\geq 1,\nonumber \\
\widetilde{q}_{n,p} & \underset{n \to \infty}{\sim}& \widetilde{C}_{p} 12^n n^{-5/2}, \label{asympqnpt}\\
\widetilde{C}_{p} & =& \frac{(3p)!}{p!(2p-1)!}3^{-p}2^{p-1} \pi^{-1/2} \underset{p \to \infty}{\sim}  \frac{\sqrt{3p}}{2\pi}\left(\frac{9}{2}\right)^p. \nonumber
\end{eqnarray}
To simplify notation we introduce $W_{c}(z)= W(12^{-1},z)$ (resp.\,$\widetilde{W}_c(z)= W(12^{-1},z))$ the generating function of quadrangulations with general (resp.\,simple) boundary taken at the critical point $g=\frac{1}{12}$ for the size.

\begin{remark} \label{universality} The $n^{-5/2}$ and $\sqrt{p}$ polynomial corrections in the asymptotics of $q_{n,p}$ and $\widetilde{q}_{n,p}$ are common features in planar structures with boundary, in particular it holds for other ``reasonable'' classes of maps with boundary such as triangulations. These exponents turn out to rule the large scale structure of such maps. 
\end{remark}

For all $n,p \geq0$, we denote by $Q_{n,p}$ and $ \widetilde{Q}_{n,p}$ random variables with uniform distributions over $ \mathcal{Q}_{n,p}$ and $ \mathcal{ \widetilde{Q}}_{n,p}$ respectively. In the next section we recall the definition of the UIPQ and construct the UIPQ with simple boundary from it.

\subsection{The UIPQ with simple boundary} \label{UIPQ:simple}

Recall the metric $ \mathrm{d_{map}}$ presented in the Introduction. The set of all finite rooted planar quadrangulations with boundary is not complete for this metric and we will have to work in its completion $ \mathcal{Q}$. The additional elements of this set are called infinite quadrangulations with boundary. Formally they can be seen as sequences $(q_1, ... , q_n, ...)$ of finite rooted quadrangulations with boundary such that for any $r\geq 0$, $B_{r}(q_{n})$ is eventually constant. See \cite{CMMinfini} for more details. Recall from \eqref{def:UIPQ} that the UIPQ is the weak limit in the sense of $ \mathrm{d_{map}}$ of uniform rooted quadrangulations whose size tends to infinity. 
\medskip

\begin{figure}[h] 
\begin{center}
\includegraphics[width=12cm]{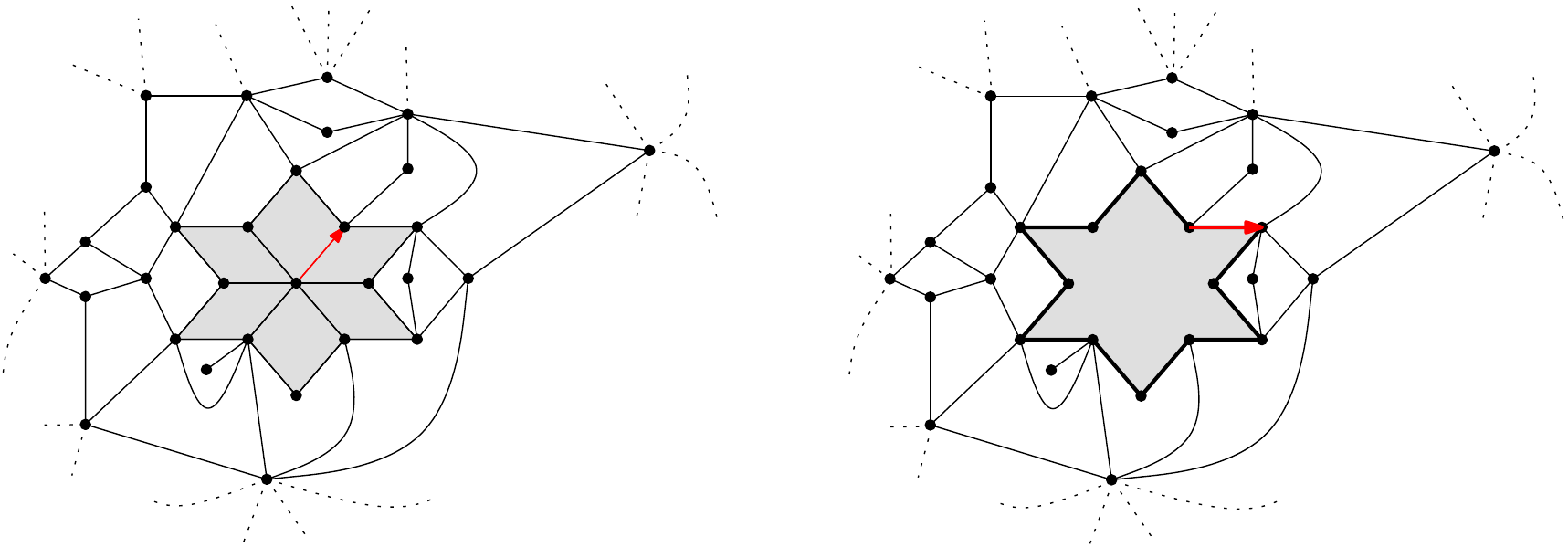} 
\caption{ \label{star}A fragment of a rooted quadrangulation with a special neighborhood of the root  and the rooted quadrangulation with simple boundary obtained by removing this neighborhood and moving the root edge.}
\end{center}
 \end{figure}

We can already use \eqref{def:UIPQ} to deduce a similar convergence result for rooted quadrangulations with a simple boundary. Indeed, notice that a rooted quadrangulation with $n$ faces and perimeter $2p$ can be turned into a rooted quadrangulation with $n+p$ faces that has a special neighborhood around the origin. More precisely, if $q_{n+p} \in \mathcal{Q}_{n+p,1}$ is a rooted quadrangulation such that the neighborhood of the root edge is composed of $p$ squares arranged like a star around the origin of the root edge as depicted in Fig.\,\ref{star} (note that the vertices on the boundary of the star must be pairwise distinct), then we can remove this star from $q_{n+p}$  and move the root edge in a deterministic way to obtain a quadrangulation with boundary of perimeter $2p$ with $n$ internal faces. This operation is reversible.

Hence  the uniform distribution over $\mathcal{ \widetilde{Q}}_{n,p}$ can be seen as the uniform distribution over $\mathcal{Q}_{n+p,1}$ conditioned on having a ``starred neighborhood'' composed of $p$ squares. For $p \geq 1$, we condition a uniform infinite planar quadrangulation $Q_{\infty}$ to have a ``starred neighborhood'' with $p$ squares (event of positive probability) and denote  by $\widetilde{Q}_{\infty,p}$ the complement of this neighborhood, which is an infinite quadrangulation with simple boundary of perimeter $2p$ rooted as explained before. 
The convergence \eqref{def:UIPQ} together with the preceding remarks then yield

 \begin{eqnarray} \label{def:UIPQs} \widetilde{Q}_{n,p} &\xrightarrow[n\to\infty]{}& \widetilde{Q}_{\infty,p},  \end{eqnarray}
 in distribution in the sense of $ \mathrm{d_{map}}$. The random variable $\widetilde{Q}_{\infty,p}$ is called the uniform infinite planar quadrangulation (UIPQ) with simple boundary of perimeter $2p$.

 \begin{remark} It is not easy to deduce from \eqref{def:UIPQ} a similar convergence result for quadrangulations with a general boundary: This is due to the fact that quadrangulations with general boundary are not \emph{rigid} in the sense of \cite[Definition 4.7]{AS03}. We prefer to take a different approach to define the UIPQ with general boundary in the next sections. \end{remark}

\section{Bijective Representation} \label{bijections}
In this section we extend the bijective approach of the UIPQ developed in \cite{CMMinfini} to the case of quadrangulations with boundary using the tools of \cite{BG09}.  For technical reasons, we will have to consider pointed quadrangulations:  A (rooted) quadrangulation with boundary is \emph{pointed} if it is given with a distinguished vertex $\rho$. We let $\mathcal{Q}_{n,p}^\bullet$ (resp.\,$\widetilde{ \mathcal{Q}}^{\bullet}_{n,p}$) be the set of all rooted pointed quadrangulations with general (resp.\,simple) boundary and $q_{n,p}^{\bullet}$ (resp.\,$\widetilde{q}^{\bullet}_{n,p}$) its cardinality. 
 
\subsection{Trees}
We use the same notation as in \cite{CMMinfini}. Let $
\mathcal{U} = \cup_{n=0}^{\infty} \mathbb{N}^n
$,
where $\mathbb{N} = \{ 1,2, \ldots \}$ and $\mathbb{N}^0 = \{
\varnothing \}$ by convention. An element $u$ of $\mathcal{U}$ is thus
a finite sequence of positive integers. If $u, v \in \mathcal{U}$,
$uv$ denotes the concatenation of $u$ and $v$. If $v$ is of the form
$uj$ with $j \in \mathbb{N}$, we say that $u$ is the \emph{parent} of
$v$ or that $v$ is a \emph{child} of $u$. More generally, if $v$ is of
the form $uw$ for $u,w \in \mathcal{U}$, we say that $u$ is an
\emph{ancestor} of $v$ or that $v$ is a \emph{descendant} of $u$. A
\emph{plane tree} $\tau$ is a (finite or infinite) subset of
$\mathcal{U}$ such that
\begin{enumerate}
\item $\varnothing \in \tau$ ($\varnothing$ is called the \emph{root}
  of $\tau$),
\item if $v \in \tau$ and $v \neq \varnothing$, the parent of $v$
  belongs to $\tau$
\item for every $u \in \mathcal{U}$ there exists $k_u(\tau) \geq 0$
  such that $uj \in \tau$ if and only if $j \leq k_u(\tau)$.
\end{enumerate}
A plane tree can be seen as a graph, in which an edge links
two vertices $u,v$ such that $u$ is the parent of $v$ or vice-versa. 
This graph is of course a tree in the graph-theoretic sense, and has a
natural embedding in the plane, in which the edges from a vertex $u$
to its children $u1,\ldots,uk_u(\tau)$ are drawn from left to right. We let $|u|$ be the length of the word $u$. The
integer $|\tau|$ denotes the number of edges of $\tau$ and is called
the size of $\tau$.  A \emph{corner} of a tree is an angular sector formed by two consecutive edges in the clockwise contour. A \emph{spine} in a tree $\tau$ is an infinite sequence
$u_{0},u_{1},u_{2},\ldots$ in $\tau$ such that $u_{0}=\varnothing$ and $u_i$
is the parent of $u_{i+1}$ for every $i\geq 0$. 
In this work, unless explicitly mentioned, all the trees considered are plane trees.

\paragraph{The uniform infinite plane tree.} For any plane tree $\tau$ and any $h \geq0$ we define the tree $[\tau]_h=\{u \in \tau : |u| \leq h\}$ as the tree $\tau$ restricted to the first $h$ generations. If $\tau$ and $\tau'$ are two plane trees, we set 
\begin{eqnarray*} \mathrm{d_{tree}}(\tau,\tau') &=& \big( 1+ \sup\{ h \geq 0 : [\tau]_h= [\tau']_h \}\big) ^{-1}. \end{eqnarray*} Obviously, $ \mathrm{d_{tree}}$ is a distance on the set of all plane trees. In the following, for every $n \geq 0$, we denote by $T_n$ a random variable uniformly distributed over the set of all rooted plane trees with $n$ edges. It is standard (see \cite{Kes86,LPP95b,Dur03}) that there exists a random infinite plane tree $T_\infty$ with one spine  called the uniform infinite plane tree, or critical geometric Galton-Watson tree conditioned to survive, such that we have the convergence in distribution for $ \mathrm{d_{tree}}$
 \begin{eqnarray}T_n & \xrightarrow[n\to\infty]{(d)} & T_\infty.  \label{def:tinfty}\end{eqnarray}
The tree $T_\infty$ can be informally described as follows. Start with a semi-infinite line of vertices (which will be the unique spine of the tree, rooted at the first vertex of the spine), then on the left and right hand side of each vertex of the spine, graft independent critical geometric Galton-Watson trees with parameter $1/2$. The resulting plane tree has the same distribution as $T_\infty$. See \cite{BCsnake,CMMinfini} for more details.
\paragraph{Labeled trees.}

A \emph{rooted labeled tree} (or spatial tree) is a pair $\theta =
(\tau, (\ell(u))_{u \in \tau})$ that consists of a plane tree $\tau$
and a collection of integer labels assigned to the vertices of $\tau$, 
such that if $u,v \in \tau$ are neighbors then $|\ell(u) - \ell(v)| \leq 1$.  Unless mentioned, the label of the root vertex is $0$. If $\theta = (\tau,\ell)$ is a labeled tree, $|\theta| = |\tau|$ is the size of $\theta$.  Obviously, the distance $ \mathrm{d_{tree}}$ can be extended to labeled trees by taking into account the labels, and we keep the notation $ \mathrm{d_{tree}}$ for this distance. 

Let $\tau$ be a random plane tree and, conditionally on $\tau$, consider a sequence of independent identically distributed random variables uniformly distributed over $ \left\{-1,0,1\right\}$ carried by each edge of $\tau$. For any vertex $u$ of $\tau$, the label of $u$ is defined as the sum of the variables carried by the edges along the unique path from the root $\varnothing$ to $u$. This labeling is called the \emph{uniform} labeling of $\tau$. When the tree $\tau$ is a geometrical critical Galton-Watson tree (conditioned to survive), we will speak of the ``uniform labeled critical geometric Galton-Watson tree (conditioned to survive)". Using the notation of \cite{CMMinfini}, we denote by $\mathscr{S}$ the set of all labeled infinite trees with only one spine such that the infimum of the labels along the spine is $-\infty$. Note that if $\tau$ is an infinite tree with one spine and $\theta= (\tau, \ell)$ is a uniform labeling of $\tau$ then $\theta \in \mathscr{S}$ almost surely.

\subsection{Treed bridges}

 A \emph{bridge} of length $2p$ is a sequence of integers $x_1, x_2, ... , x_{2p}$ such that $x_1=0$ and for every $i \in \{1, ... 2p\}$ we have $|x_{i}-x_{i+1}|=1$, where by convention we let $x_{2p+1}=x_1$. 
Note that in any bridge of length $2p$, there are exactly $p$ \emph{down-steps}, which are the indices $i\in \{1,2,\ldots,2p\}$ such that $x_{i+1}=x_i-1$. A \emph{labeled treed bridge} of  size $n$ and length $2p$ is a bridge $\mathrm{b}_p=(x_1, ... , x_{2p})$ together with $p$  non-empty labeled plane trees $\theta_1, ... , \theta_p$ (with root label $0$) such that the sum of the sizes of the trees $\theta_1, ... , \theta_p$ is $n$. We denote by $\mathscr{B}_p$ the union set of all labeled finite treed bridges and infinite labeled treed bridges $(\mathrm{b}_p; \theta_1, ... , \theta_p)$ such that \emph{one and only one} of the labeled trees $\theta_i$ belongs to $\mathscr{S}$, all  others being finite. In the following, unless explicitly mentioned, all labeled treed bridges considered  belong to $\mathscr{B}_p$ for some $p\geq 1$.

\paragraph{Representation.}  Let $ \mathbf{b}=(\mathrm{b}_p; \theta_1, ... , \theta_p)$ be a treed bridge of $\mathscr{B}_p$. If $\mathrm{b}_p=(x_1,... , x_{2p})$ we denote $(x_{i_1},x_{i_1+1}), ... , (x_{i_p},x_{i_p+1})$ its $p$ down-steps. We  construct \emph{a representation} of $\mathbf{b}$ in the plane as follows. Let $\mathcal{C}$ be a proper embedding in the plane of a cycle of length $2p$. We label the vertices of $\mathcal{C}$ starting from a distinguished vertex in the clockwise order by the values of the bridge $\mathrm{b}_p$. Now we graft (proper embeddings of) the trees $\theta_1, ..., \theta_p$ in the infinite component of $\mathbb{R}^2 \backslash \mathcal{C}$ in such a way that the tree $\theta_k$ is grafted on the $i_k$th point of $\mathcal{C}$ corresponding to the value $x_{i_k}$ and we shift all the labels of this tree by $x_{i_k}$. This representation can be constructed in such a way that the embedding is proper (no edges are crossing except possibly at their endpoints) and such that the sequence of vertices of the embedding has no accumulation point in $\mathbb{R}^2$ (recall that there is at most one infinite tree with only one spine). See Fig.\,\ref{representation} below.  \begin{figure}[h] \begin{center}
    \includegraphics[width=14cm]{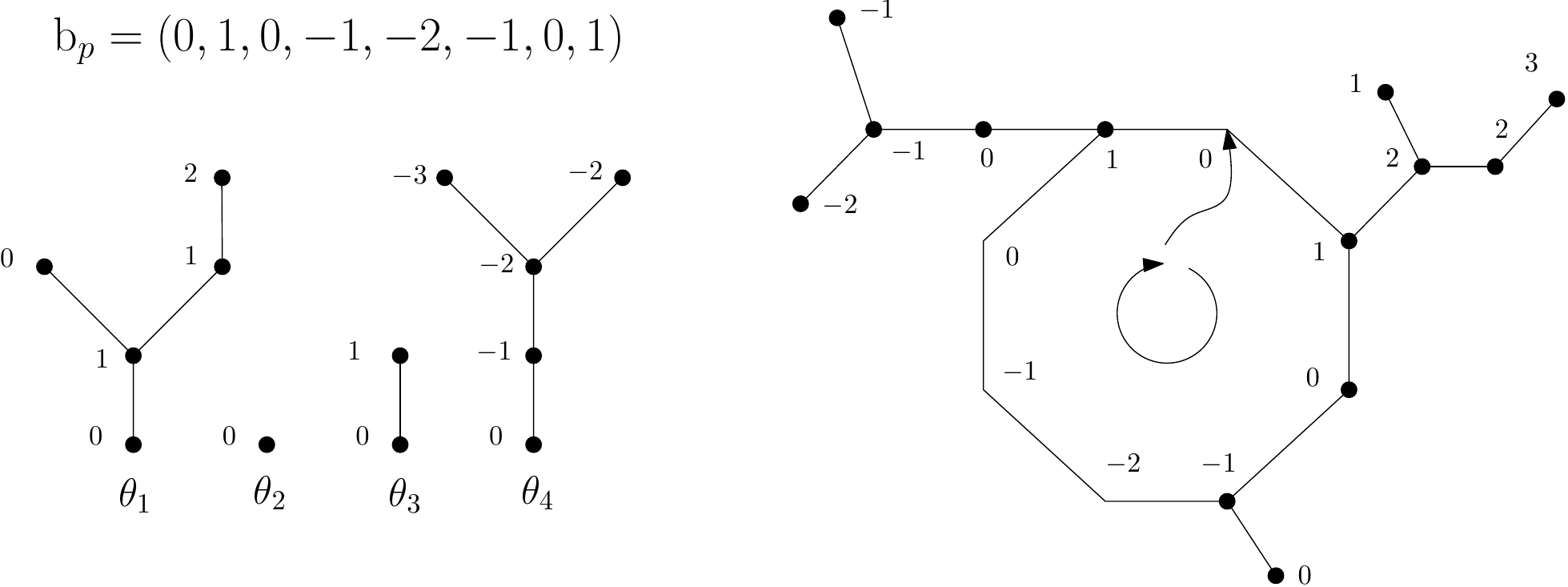} \\
    \caption{ \label{representation} A labeled finite treed bridge and an representation of it. The black dots represent the vertices of the trees of the treed bridge.}  \end{center}
 \end{figure}

 The vertex set of this representation is thus formed by the union of the vertices of the trees $\theta_1, ... , \theta_p$ and of the vertices of the cycle which are not down-steps. The labeling of these vertices, which is given by the bridge on the cycle and the shifted labelings of the trees is denoted by $\ell_\mathbf{b}$. We will often abuse notation and write $u \in \theta_i$ for a vertex in the representation that belongs to the embedding of the tree $\theta_i$.\\
 Recall that a corner of a proper embedding of a graph in the plane is an angular sector formed by two consecutive edges in clockwise order. In the case of a representation of a labeled treed bridge we can consider the set of corners of the infinite component of the plane minus the embedding. This set, although possibly infinite, inherits a cyclic order from the clockwise order of the plane.  The label of a corner is that of its attached vertex.

\paragraph{The uniform infinite labeled treed bridge.}
Let $p\geq 1$. We say that a sequence of  labeled treed bridges $(\mathrm{b}_p^{(n)}; \theta_1^{(n)}, ... , \theta_p^{(n)})$ of length $2p$ converges to $(\mathrm{b}_p; \theta_1,... , \theta_p)$, if eventually  $\mathrm{b}_p^{(n)}=\mathrm{b}_p$ and $\theta_i^{(n)}$ converges towards $\theta_i$ for any $i \in \{1, ... ,p\}$ with respect to $ \mathrm{d_{tree}}$. This convergence is obviously metrizable by the metric 
 \begin{eqnarray*} \mathrm{d_{bridge}}\big((\mathrm{b}_p; \theta_1, ... , \theta_p),(\mathrm{b}'_p; \theta'_1, ... , \theta'_p)\big) &=& \mathbf{1}_{\mathrm{b}_p=\mathrm{b}'_p} \sup_{1 \leq i \leq p} \mathrm{d_{tree}}(\theta_i, \theta'_i) + \mathbf{1}_{ \mathrm{b}_p \ne \mathrm{b}_{p'}}.  \end{eqnarray*}
In the rest of this work, $\mathbf{B}_{n,p} = (B_p; \Theta_1^{(n)}, ... , \Theta_p^{(n)})$ is a uniform labeled bridge of size $n$ and length $2p$. Note that  $B_p$ is uniformly distributed over the set of bridges of length $2p$. We introduce the analog of the labeled critical geometric Galton-Watson tree conditioned to survive in the setting of treed bridges. The \emph{uniform infinite labeled treed bridge} denoted by $\mathbf{B}_{\infty,p}=(B_p, \Theta_1, ... , \Theta_p)$ is constructed as follows. Firstly, $B_p$ is a uniform bridge of length $2p$. Then choose $i_{0} \in \{1,... ,p\}$ uniformly and independently of $B_p$. Conditionally on $B_p$ and $i_{0}$,  the trees $\Theta_1, ... , \Theta_p$ are independent, $\Theta_{i_{0}}$ being  a uniform labeled critical geometric Galton-Watson tree conditioned to survive and all other  $\Theta_j$ for $j\ne i$ are uniform labeled critical geometric Galton-Watson trees. Notice that $\mathbf{B}_{\infty,p}$ almost surely belongs to $\mathscr{B}_p$. Then the analog of \eqref{def:tinfty} becomes:

\begin{proposition}\label{Kestenbis} 
We have the following convergence in distribution for $ \mathrm{d_{bridge}}$
\begin{eqnarray}
\mathbf{B}_{n,p} & \xrightarrow[n\to\infty]{(d)} & \mathbf{B}_{\infty,p}.
\end{eqnarray}
\end{proposition}
\proof Since conditionally on the structure of the trees, the bridge and the labels are uniform they do not play any crucial role in the convergence. We just have to prove that if $\tau_1^{(n)}, ... , \tau_p^{(n)}$ are $p$ plane trees chosen uniformly among all $p$-uples of  plane trees such that $|\tau_1^{(n)}|+ ... + |\tau_p^{(n)}|=n$ then we have the following weak convergence 
\begin{eqnarray*} (\tau_1^{(n)}, ... , \tau_p^{(n)}) & \xrightarrow[n\to\infty]{}& (\tau_1, ... , \tau_p), \end{eqnarray*} where the  distribution of $(\tau_1, ... , \tau_p)$ is described as follows: Choose $i_{0}$ uniformly among $\{1, ... , p\}$, then conditionally on $i_{0}$ the $\tau_k$'s are independent, $\tau_{i_0}=T_{\infty}$ in distribution and the other trees are critical geometric Galton-Watson trees. This fact is 
standard but we include a proof for the reader's convenience. For $p \geq 1$ and $n\geq 0$, we let $\mathrm{Cat}(n,p)$ be the number of finite sequences (``forests'') of $p$ trees with $n$ edges in total, so that by a well-known formula (see e.g.\ \cite{Pit06})
 \begin{eqnarray} \mathrm{Cat}(n,p) \quad = \quad
\frac{p}{2n+p}\binom{2n+p}{n}
\quad  \underset{n \to \infty}{\sim} \quad  4^n n^{-3/2}\frac{2^{p-1}p}{\sqrt{\pi}} \label{asympcatnp}.
\end{eqnarray}
We also let $\mathrm{Cat}(n)=\mathrm{Cat}(n,1)$. 
Let $f_1, ... , f_p$ be bounded continuous functions for $ \mathrm{d_{tree}}$. By definition of the distribution of $(\tau_1^{(n)},... , \tau_p^{(n)})$ we have 
\begin{eqnarray}
{E}\left[\prod_{i=1}^p f_i(\tau_i^{(n)})\right] & = & \frac{1}{ \mathrm{Cat}(n,p)} \sum_{n_1+ ... + n_p=n} \prod_{i=1}^p \mathrm{Cat}(n_i) {E}[f_i(T_{n_i})], \label{exacttree}
\end{eqnarray} 
where $T_{n_i}$ denotes a uniform plane tree on $n_i$ edges.  We first estimate the probability that two of the trees $\tau_1^{(n)}, ... , \tau_p^{(n)}$ have a size larger than some large constant $a>0$. For that purpose we recall a classical lemma whose proof is very similar to \cite[Lemma 2.5]{AS03} and is left to the reader.\begin{lemma} \label{utile} For any $ \beta>1$, and any $p \geq 0$ there exists a constant $c(\beta,p)$ such that for any $a>0$ and any $n\geq 0$ we have
$$\begin{array}{rcl}

 \displaystyle \sum_{\begin{subarray}{c}
n_1 + ... + n_p = n \\
n_1,n_2>a
\end{subarray}
} \left(\prod_{i=1}^p (n_i+1)^{-\beta}\right)& \leq& c(\beta,p) n^{-\beta}a^{-\beta+1} .\\

\end{array}$$
\end{lemma}
Using the asymptotic \eqref{asympcatnp} with $p=1$, we deduce that there exists a constant $c_1$ such that for any $n \geq 1$ we have $\mathrm{Cat}(n) \leq c_1 4^n(n+1)^{-3/2}$. We can thus use Lemma \ref{utile} to deduce that there exists a constant $c_2>0$ such that for every $n\geq 1$, the probability that two of the trees $\tau_1^{(n)},..., \tau_p^{(n)}$ have size larger than $a>0$ is less than $ c_2a^{-1/2}$. Hence, for large $n$'s the right-hand side of \eqref{exacttree} becomes 
$$ \varepsilon_{n,a}+ \sum_{i=1}^p \left( \sum_{0 \leq n_{1}, ..., \widehat{n_{i}}, ..., n_{p} \leq a}\frac{\mathrm{Cat}\left(n-\sum_{j\ne i} n_j\right)}{ \mathrm{Cat}(n,p)}{E}[f_{i}(T_{n-\sum_{j\ne i }n_j})] \prod_{j \ne i} \mathrm{Cat}(n_j) {E}[f_j(T_{n_j})]    \right),$$ where $ \varepsilon_{n,a} \leq c_2 a^{-1/2}$ uniformly in $n \geq 1$. Moreover, \eqref{def:tinfty}  implies that ${E}[f(T_n)] \to {E}[f(T_\infty)]$ for any bounded continuous functional for $\mathrm{d_{tree}}$ as $n \to \infty$. So, using once more the asymptotic \eqref{asympcatnp}, we can let $n \to \infty$  followed by $a \to \infty$ in the last display and obtain 
\begin{eqnarray}
{E}\left[\prod_{i=1}^p f_i(\tau_i^{(n)})\right] & \xrightarrow[n\to\infty]{} & \sum_{i=1}^p \frac{1}{p}{E}[f_i(T_{\infty})] \prod_{j \ne i} \sum_{n_j =0}^\infty\frac{\mathrm{Cat}(n_j)}{2\cdot 4^{n_j}} {E}[f_j(T_{n_j})].
\end{eqnarray}
The sum over indices $n_j$ is $E[f_j(T)]$, where $T$ is a uniform  critical geometric Galton-Watson tree, which proves the desired result.\endproof

\subsection{From treed bridges to quadrangulations with boundary}
\subsubsection{Finite bijection}
The bijection presented in this section is taken from \cite{BG09} with minor adaptations. This is a one-to-one correspondence between, on the one hand, the set $\mathcal{Q}_{n,p}^\bullet$ of all 
rooted and pointed quadrangulations with boundary of perimeter $2p$ and size $n$ and, on the other hand, the set of all  labeled treed bridges of length $2p$ and size $n$. 
We only present the mapping from  labeled treed bridges to quadrangulations, the reverse direction can be found in \cite{BG09}. \\

\begin{figure}[h]
 \begin{center}
\includegraphics[width=12cm]{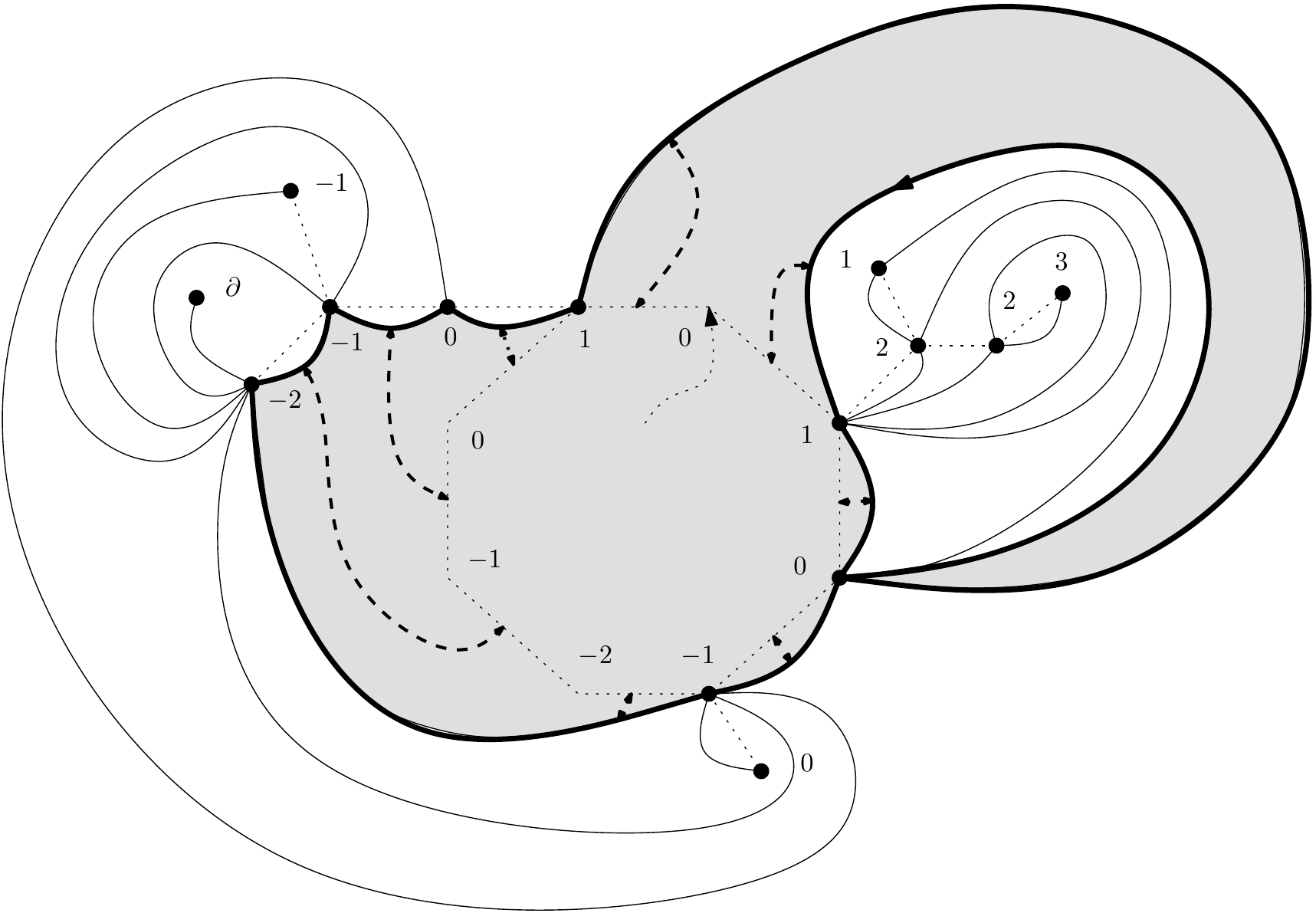} \\
\caption{ \label{BG} Construction of the rooted and pointed quadrangulation with boundary associated with a treed bridge. Note the correspondence between the embedded cycle and the boundary of the map.}
\end{center}
\end{figure}

Let $\mathbf{b}=(\mathrm{b}_p; \theta_1, ... , \theta_p)$ be a  labeled treed bridge of size $n$ and perimeter $2p$. We consider a representation $\mathcal{E}$ of this treed bridge in the plane. Recall that the labeling of the vertices of this representation is denoted by $\ell_\mathbf{b}$. 
Let $\mathscr{C}$ be the set of corners of the infinite component of $\mathbb{R}^2 \backslash \mathcal{E}$ that are also incident to vertices belonging of the grafted trees, that is, we erase the corners coming from angular sectors around the vertices of the cycle that are not down steps of the bridge (see Fig.\,\ref{BG}). This set inherits a clockwise cyclic order. 

We now associate a quadrangulation with $ \mathbf{b}$ by the following device: We start by putting an extra vertex denoted $\rho$ in the infinite component of $ \mathbb{R}^2 \setminus \mathcal{E}$. Then for each corner $c \in \mathscr{C}$, we draw an edge between $c$ and the first corner $c'\in \mathscr{C}$ for the clockwise order such that the label $\ell_\mathbf{b}(c')$ of $c'$ equals $\ell_\mathbf{b}(c)-1$: This corner is called \emph{the successor} of $c$. If there is no such corner (this happens only if the corner $c$ has minimal label) then we draw an edge between $c$ and $\rho$. This construction can be done in such a way that the edges are non-crossing. After erasing  the representation $\mathcal{E}$ of $\mathbf{b}$ we obtain a quadrangulation $q$ with size $n$ and a boundary whose vertex set is the union of the vertices of (the embeddings of) $\theta_i$ for $i \in \{1, ... , p\}$ plus the extra vertex $\rho$.

Note that there is a one-to-one order-preserving correspondence between the edges of the cycle of $\mathcal{E}$ and the edges of the external face of $q$, see Fig.\,\ref{BG}. The distinguished oriented edge $\vec{e}$ of $q$ is the edge that corresponds to the first step of the bridge oriented such that the external face is on the right-hand side of $\vec{e}$.  We denote the rooted  quadrangulation $(q,\vec{e})$ pointed at $\rho$ by $\Phi(\mathbf{b})$.  
Furthermore, using the identification of the vertices of the map $\Phi( \mathbf{b})$ with $\cup_{i }\theta_{i} \cup \{\rho\}$, for every $u \in \cup_{i} \theta_{i}$ we have 
\begin{eqnarray} \ell_{ \mathbf{b}}(x)-\ell_{ \mathbf{b}}(\rho) &=& \mathrm{d_{gr}}(x,\rho), \label{interpret}\end{eqnarray}
where $ \mathrm{d_{gr}}$ is the graph distance in $ \Phi( \mathbf{b})$.
\label{unconditioned}
\subsubsection{Extended construction}
We can extend the preceding mapping $\Phi$ 
to the case when the treed bridge $\mathbf{b}$ is infinite but still belongs to $\mathscr{B}_p$. The extension is very similar to that of \cite{CMMinfini}. Basically, the construction goes through.  The only point that is changed is that every corner attached to a tree in the infinite component of the embedding will find a successor, that is, there is no  need to add an extra vertex $\rho$. The extended mapping that we  denote $ \Phi$, associates with every labeled treed bridge in $\mathscr{B}_p$ an infinite rooted quadrangulation $q$ with boundary of perimeter $2p$ (this quadrangulation is not pointed anymore) whose vertex set is the union of the vertices of the trees of the bridge. The correspondence between the cycle of the representation and the boundary of the map is still preserved. The shifted labels lose their finite interpretation \eqref{interpret} (see \eqref{interpret2} in Theorem \ref{CMM10bound}).  However, since any neighboring vertices in the quadrangulation $q$  have labels that differ in absolute value by exactly $1$, we deduce that for any vertices $u,v$ in the resulting quadrangulation we have (with the identification of the vertices of the quadrangulation with those of the trees of $\mathbf{b}$)
\begin{eqnarray} \mathrm{d}_{\mathrm{gr}}(u,v) &\geq& |\ell_\mathbf{b}(u)-\ell_\mathbf{b}(v)|.  \label{bound} \end{eqnarray}
\begin{proposition} \label{continuity}
The extended Schaeffer mapping $\Phi :\mathscr{B}_p \longrightarrow \mathcal{Q}$ is continuous with respect to the metrics $ \mathrm{d_{bridge}}$ and $ \mathrm{d_{map}}$.
\end{proposition}
\proof The proof is similar to that of \cite[Proposition 1]{CMMinfini} and is left to the reader.\endproof

\section{The UIPQ with general boundary} \label{sec:UIPQbound}
\subsection{Construction}
The following theorem in an extension of the main result of \cite{CMMinfini} to the case of quadrangulations with boundary. Recall that $Q_{n,p}$ is uniformly distributed over $\mathcal{Q}_{n,p}$.

\begin{theorem} \label{CMM10bound} $(i)$\ For any $p \in \{ 1,2,3,...\}$ we have the following convergence in distribution for $ \mathrm{d_{map}}$
\begin{eqnarray}
Q_{n,p}& \xrightarrow[n\to\infty]{(d)}& Q_{\infty,p}, \nonumber 
\end{eqnarray}
where $Q_{\infty,p}$ is called \emph{the uniform infinite planar quadrangulation with boundary of perimeter $2p$}. If $ \mathbf{B}_{\infty,p}$ is a uniform infinite labeled treed bridge of length $2p$ then $Q_{\infty,p}= \Phi( \mathbf{B}_{\infty,p})$ in distribution. 

$(ii)$\ If $Q_{\infty,p}= \Phi( \mathbf{B}_{\infty,p})$ then, with the identification of the vertices of $Q_{\infty,p}$ with those of the trees of $ \mathbf{B}_{\infty,p}$, we have for any $u,v \in Q_{\infty,p}$
\begin{eqnarray} \lim_{z \to \infty} \big(\mathrm{d_{gr}}(u,z)-\mathrm{d_{gr}}(v,z)\big) &=& \ell_{\mathbf{B}_{\infty,p}}(u)-\ell_{\mathbf{B}_{\infty,p}}(v). \label{interpret2}\end{eqnarray}
\end{theorem}

\proof[First part of Theorem \ref{CMM10bound}] An application of Euler's 
formula shows that every quadrangulation with a boundary of perimeter $2p$ and size $n$ has exactly $n+p+1$ vertices. We deduce from the preceding section that after forgetting the distinguished point, the rooted quadrangulation $\Phi(\mathbf{B}_{n,p})$ is uniform over $\mathcal{Q}_{n,p}$. The first part of the theorem then follows from Proposition \ref{Kestenbis} and Proposition \ref{continuity}. \\
For the second part of the theorem, notice that the case $p=1$ is proved in \cite{CMMinfini}. The general case will follow from this case using some surgical operation that we present in Section \ref{sec:pruning}.\endproof
\begin{remark} In the construction $Q_{\infty,p}= \Phi( \mathbf{B}_{\infty,p})$, since the origin of the root edge of $Q_{\infty}$ automatically has label $0$, the formula \eqref{interpret2}  can be used  to recover $ \ell_{ \mathbf{B}_{\infty,p}}$ as a measurable function of $ Q_{\infty,p}$. Using an extension of the reversed construction ``$\Phi^{-1}$'' (see \cite{BG09}) it can be proved following the lines of \cite{CMMinfini} that the treed bridge $ \mathbf{B}_{\infty,p}$ itself can be recovered as a measurable function of $Q_{\infty,p}$. We leave this to the interested reader.
\end{remark}

\begin{remark} Recall that the UIPQ with a simple boundary was defined from the standard UIPQ in Section \ref{UIPQ:simple}. It is also possible to define the UIPQ with simple boundary $\widetilde{Q}_{\infty,p}$ as the variable $Q_{\infty,p}$ conditioned on having a simple boundary. However, thanks to the asymptotics \eqref{asympqnp} and \eqref{asympqnpt} the probability that the boundary of $Q_{\infty,p}$ is simple is easily seen to be 
$$ {P}(Q_{\infty,p} \mbox{ is simple}) \quad =\quad \frac{ \widetilde{C}_p}{ C_p} \quad \underset{p \to \infty}{\sim} \quad  \sqrt{3}\left(\frac{9}{16}\right)^p.$$
This exponential decay is not useful if we want to derive properties of $\widetilde{Q}_{\infty,p}$ from properties of $Q_{\infty,p}$ for large $p$'s. For that purpose we develop in Section \ref{sec:pruning} another link between $\widetilde{Q}_{\infty,p}$ and ${Q}_{\infty,p}$ based on a pruning procedure.
\end{remark}
\subsection{Aperture of the UIPQ with general boundary}
Recall that if $q$ is a quadrangulation with boundary, the \emph{aperture} of $q$ is the maximal graph distance between any two points of the boundary 
 \begin{eqnarray*} \mathrm{aper}(q) &=& \max \{ \mathrm{d_{gr}}(u,v): u,v \in \partial q\}. \end{eqnarray*}
The main result of this section provides bounds on the  aperture of $Q_{\infty,p}$ for large $p$'s. It is based on the construction $Q_{\infty,p} = \Phi( \mathbf{B}_{\infty,p})$ of Theorem \ref{CMM10bound} and on properties of some specific geodesics in $Q_{\infty,p}$.

\begin{figure}[!h]
 \begin{center}
 \includegraphics[width=12cm]{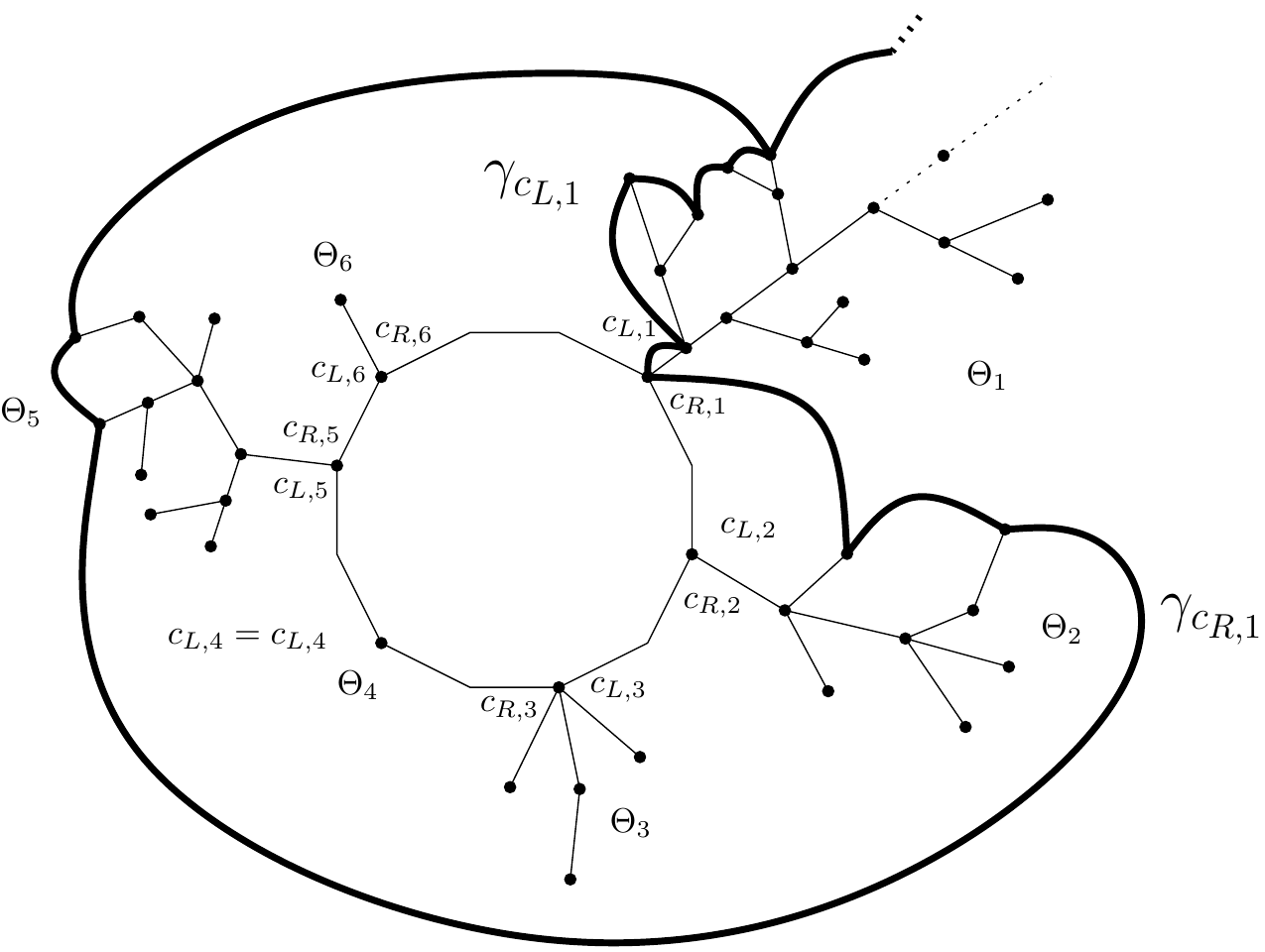}
 \caption{ \label{corners} Left and right most corners of the trees $\Theta_1, ... , \Theta_6$ and the geodesics $\gamma_{c_{L,1}}$ and $\gamma_{c_{R,1}}$.}
 \end{center}
 \end{figure}

\begin{theorem} \label{diamext} The  aperture of $Q_{\infty,p}$ is exponentially concentrated around the order of magnitude $\sqrt{p}$. More precisely, there exist $c_1,c_2>0$ such that for all $\lambda >0$ and every $p \in \{1,2,3,... \}$ we have
\begin{eqnarray*} {P}\big( \mathrm{aper}(Q_{\infty,p}) \notin [\lambda^{-1} \sqrt{p}, \lambda \sqrt{p}] \big)  &\leq& c_1 \exp(-c_2\lambda^{2/3}). \end{eqnarray*}
\end{theorem}
Fix $p \in \{1,2,3, ...\}$.  To simplify notation, we write $\mathbf{B}=(B_p; \Theta_1, ... , \Theta_p) \in \mathscr{B}_p$ for the uniform infinite labeled treed bridge of length $2p$ and assume that $Q_{\infty,p}= \Phi(\mathbf{B})$. We will always use a representation $\mathcal{E}$ of $\mathbf{B}$ and identify the trees $\Theta_1, ... , \Theta_p$ and the bridge $B_p$ with their embeddings in the representation, see Fig.\,\ref{representation}. Recall that the shifted labels of the trees are denoted by $\ell_{\mathbf{B}}$. We denote by $\mathscr{C}$ the set of corners of $\mathcal{E}$ that are associated with some vertex of $\cup_{i} \Theta_i$. If $c,c' \in \mathscr{C}$ we denote by $[c,c']$ the set of corners of $\mathscr{C}$ that are in-between $c$ and $c'$ for the clockwise order. For $1 \leq i \leq p$, the set of corners attached to the tree $\Theta_i$ is  $[c_{L,i},c_{R,i}]$ where $c_{L_i}$ and $c_{R,i}$ respectively denote the left and right most corner of the root vertex of $\Theta_i$ in $\mathcal{E}$, see Fig.\,\ref{corners}.

\paragraph{Simple geodesic.} We recall the notion of simple (or maximal) geodesic, see \cite[Definition~3]{CMMinfini}.  Let $c\in \mathscr{C}$. We can construct a path in the quadrangulation $\Phi(\mathbf{B})$  by starting with the corner $c$ and following iteratively its successors. This path is called \emph{the simple geodesic starting from $c$} and is denoted by $\gamma_{c}$, see Fig.\,\ref{corners} for examples. It is easy to see, thanks to \eqref{bound}, that this path is actually a geodesic in the quadrangulation. If $ \mathbf{B}$ were finite, this path would eventually end at $\rho$. In general, if $c,c'\in \mathscr{C}$, then $\gamma_c$ and $\gamma_{c'}$ merge at a corner $c''$  of label 
$$ \ell_\mathbf{B}(c'') = \max  \left \{\min_{[c,c']} \ell_\mathbf{B}, \min_{[c',c]} \ell_\mathbf{B} \right \}-1.$$

\proof[Proof of Theorem \ref{diamext}] 
We will suppose, without loss of generality, that $\Theta_1$ is the uniform labeled infinite tree so that $\Theta_2, ... , \Theta_p$ are uniform labeled critical geometric Galton-Watson trees. Let us start with a preliminary observation.\\

\textsc{Warmup.} Imagine that we construct the two simple geodesics $\gamma_{c_{L,1}}$ and $\gamma_{c_{R,1}}$ starting from the extreme corners of the root of $\Theta_1$ and denote by $\mathscr{C}$ the cycle they form until their meeting point, see Fig.\,\ref{corners}. We let $M= \max\{ \ell_{\mathbf{B}}(u) : u \in \cup_{i=2}^p\Theta_i  \mbox{ or } u = \mbox{root}(\Theta_1)\}$ and $m= \min\{ \ell_{\mathbf{B}}(u) : u \in \cup_{i=2}^p\Theta_i  \mbox{ or }u = \mbox{root}(\Theta_1)\}$. By the remark made on simple geodesics, one sees that the length of $\mathscr{C} $ is less than $2(M-m +1)$ and that every vertex in the external face of $Q_{\infty,p}$ is linked to $\mathscr{C}$ by a simple geodesic of length less than $M-m+1$. Thus we deduce that 
\begin{eqnarray} \mathrm{aper}(Q_{\infty,p}) &\leq& 3(M-m+1). \label{ineq1}\end{eqnarray} Although $M$ and $m$ are typically of order $\sqrt{p}$, yet it is possible that with a probability of order $p^{-2}$, a specific tree, say $\Theta_2$, has a height larger than $p^2$ and thus contains labels of order $\pm p$. If that happens then $M-m$ becomes of order $p$ and not $\sqrt{p}$ anymore. Thus the exponential concentration presented in Theorem \ref{diamext} cannot follow from \eqref{ineq1}.  The idea is to modify the cycle $\mathscr{C}$ in order to bypass the large trees among $\Theta_2, ... , \Theta_p$.

\textsc{Bridge.} Since $Q_{\infty,p}$ is constructed from ${\mathbf{B}}$,  we know that the edges of its boundary are in correspondence with the edges of the cycle of $\mathcal{E}$, in particular the $\ell_{ \mathbf{B}}$-labeling of the vertices of $\partial Q_{\infty,p}$ corresponds to the values $(X_1, ... , X_{2p})$ of the bridge $B_p$. From the lower bound \eqref{bound} we deduce that if $\Delta_p = \max X_i -\min X_i$ we have
\begin{eqnarray} \mathrm{aper}(Q_{\infty,p}) &\geq& \Delta_p.  \label{lowerbound}\end{eqnarray}
Since $B_p$ is a uniform bridge with $2p$ steps, classical results (which easily follow from the arguments of \cite{KM09}) show that there exist positive constants $c_1$ and $c_2$ such that for all $p \geq 0$ we have 
\begin{eqnarray} {P}\big(\Delta_p \notin [\lambda^{-1} \sqrt{p}, \lambda \sqrt{p}] \big) &\leq&  c_1\exp(-c_2 \lambda^2). \label{expbridge}\end{eqnarray} 

\textsc{Shortcut.} As we said, the strategy is now to build a simple path surrounding the external face in a very similar fashion as in the warmup but to shortcut large trees. Let us be precise.
For every $i \in \{1, ... ,p\}$, declare the tree $\Theta_i$ ``good" if the maximal displacement of the labels in $\Theta_i$ is in absolute value less than $\sqrt{p}$, that is if  \begin{eqnarray*}\sup_{u,v \in \Theta_i} |\ell_{\mathbf{B}}(u)-\ell_{\mathbf{B}}(v)| \leq \sqrt{p}. \end{eqnarray*} Call $\Theta_i$ ``bad" otherwise. For any $i \in \{2, ... , p\}$,  the probability that $\Theta_i$ is bad is less than $c_3 p^{-1}$, for some $c_3>0$ (see e.g.\,\cite[Lemma 12]{CMMinfini}). Hence, if $K$ is the number of bad trees among $\{\Theta_2, ... , \Theta_p\}$, then for any $1 \leq k \leq p$  we have 
\begin{eqnarray} {P}(K \geq k) \quad\leq \quad{p \choose k} \left(\frac{c_3}{p} \right)^k \quad\leq \quad\frac{c_3^k}{k!}. \label{K}\end{eqnarray}
We now construct the path shortcutting these trees. Recall that for $1 \leq i \leq p$, we denote $c_{L,i},c_{R,i}$ respectively the left-most and right-most corners of the root vertex of $\Theta_i$ in $\mathcal{E}$. We start with $\gamma_{c_{R,1}}$ and move along it. As soon as $\gamma_{c_{R,1}}$ meets a bad tree $\Theta_i$, we proceed as follows. From $c_{L,i}$ we start the simple geodesic $\gamma_{c_{L,i}}$. We know that it requires less than $\sqrt{p}+ \Delta_p +1$ steps for $\gamma_{c_{L,i}}$ to merge with $\gamma_{c_{R,1}}$ (which happens in $\Theta_i$), then we bypass $\Theta_i$ by considering the path formed of the beginning of $\gamma_{c_{R,1}}$ until it reaches $\Theta_i$ then go backwards along $\gamma_{c_{L,i}}$ to reach the root of $\Theta_i$ and finally continue the process with $\gamma_{c_{R,i}}$. Since $\Theta_1$ is obviously bad we also shortcut it. See Fig.\,\ref{shortcut} below.  \begin{figure}[h] \begin{center}
    \includegraphics[width=9cm]{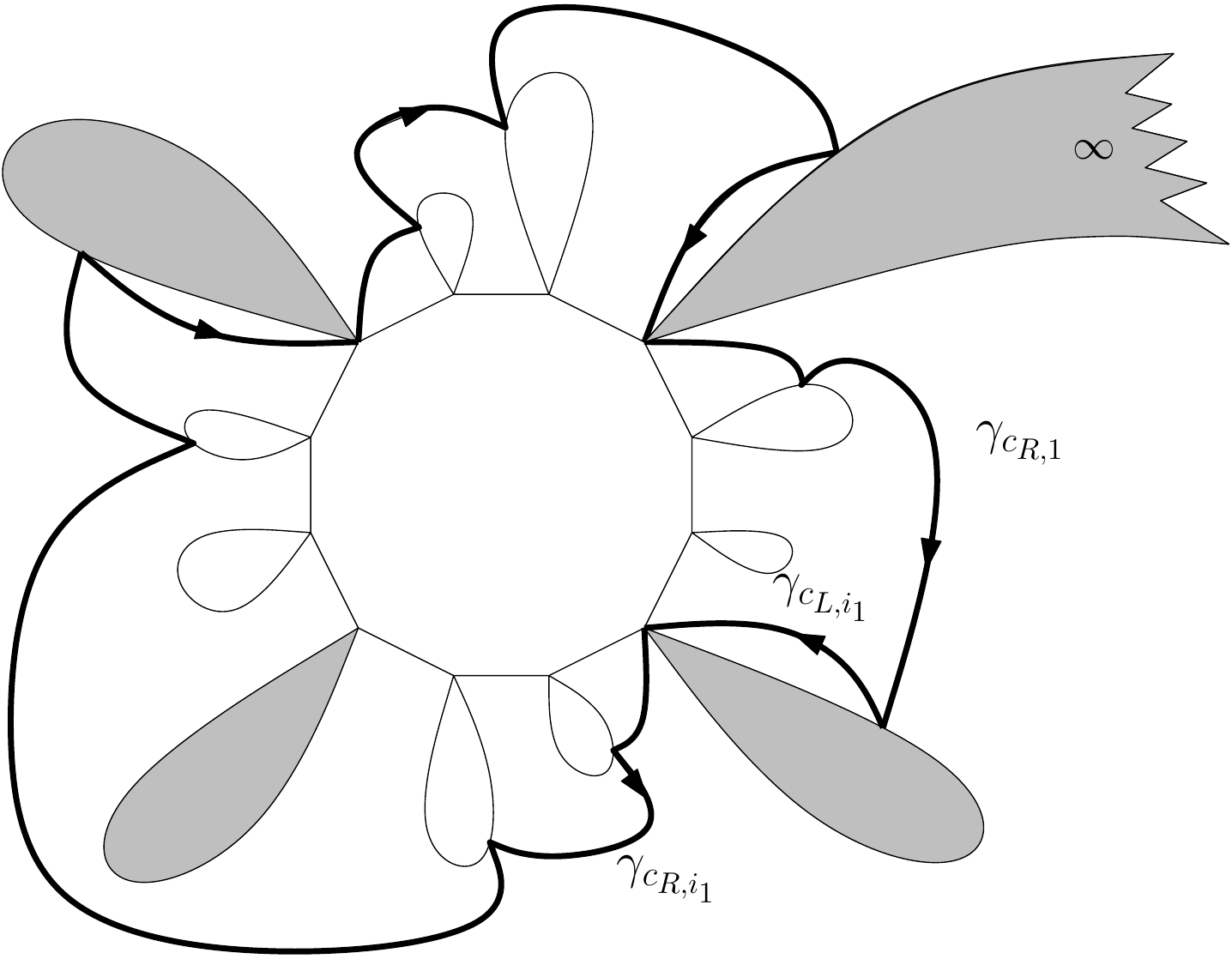}\\
    \caption{ \label{shortcut} An illustration of the proof, the bad trees are represented in gray and the cycle constructed out of simple geodesics shortcutting them in fat black line.}  \end{center} \end{figure}

At the end of the process we get a simple cycle denoted by $\mathscr{C}_{\mathrm{short}}$, which
surrounds the external face of $Q_{\infty,p}$ and whose length is at most 
\begin{eqnarray} \mathrm{Length}(\mathscr{C}_{\mathrm{short}}) &\leq& 2(K+1) (\sqrt{p}+ \Delta_p+1). \end{eqnarray}Furthermore, similarly as in the warmup part, it is easy to see that every vertex of the external face of $Q_{\infty,p}$ can be connected to $\mathscr{C}_{\mathrm{short}}$ via a simple geodesic of length less than $\sqrt{p} + \Delta_p +1$, thus the aperture of $Q_{\infty,p}$ is less than $(2K+4)\cdot(\sqrt{p} + \Delta_p +1)$, which together with \eqref{lowerbound} gives 
\begin{eqnarray}\begin{array}{ccccc}
\Delta_p &\leq &\mathrm{aper} (Q_{\infty,p}) &\leq &(2K+4)  (\sqrt{p}+\Delta_p+1).
\end{array}
\end{eqnarray}
Now, for $\lambda >1$ and $p\geq 1$, we get from the previous display   that \begin{eqnarray*}{P}(\mathrm{aper}(Q_{\infty,p}) \notin [\lambda^{-1}\sqrt{p}, \lambda \sqrt{p}]) &\leq& {P}(\Delta_p \notin [\lambda^{-1/3} \sqrt{p}, \lambda^{1/3}\sqrt{p}]) + {P}(K \geq \lambda^{2/3}/20). \end{eqnarray*} We use \eqref{expbridge} and \eqref{K} to see that the probabilities of the right hand side are of bounded by $c_4 \exp(-c_5 \lambda^{2/3})$ for some constants $c_4,c_5 >0$, which completes the proof of the theorem.\endproof

\section{Pruning} \label{sec:pruning}

\subsection{Pruning of $Q_{\infty,p}$}
Recall from Section \ref{warmup} that we can decompose a quadrangulation with a general boundary into the irreducible component containing the root edge on which quadrangulations with general boundary are attached. We now aim at a decomposition with respect to a ``big" irreducible component, which is not necessarily the one designated by the root edge (since the root edge can be located on a small component). We call this operation \emph{pruning}.

Let $q$ be a rooted quadrangulation with a boundary. Suppose that there is a unique largest irreducible component in $q$. We call this irreducible component the \emph{core} of $q$ and denote it by $ \mathrm{Core}(q)$. Attached to this core we find quadrangulations with general boundary denoted by $ \mathrm{Part}_1(q), \mathrm{Part}_2(q), ... , \mathrm{Part}_{2p}(q)$ in clockwise order where $2p$ is the perimeter of the core. Note that some of these components can be reduced to the vertex map $\dag$.  The components attached to the core are rooted at their last oriented edge visited during a clockwise contour of the external face (keeping the external face on the right). The core is either rooted at the original root edge of the map if it lies on the boundary of the core, or on the oriented edge preceding the component carrying the root edge (in that case this component is not empty). See Fig. \ref{pruningfig} below.  \begin{figure}[h] \begin{center} \includegraphics[width=16cm]{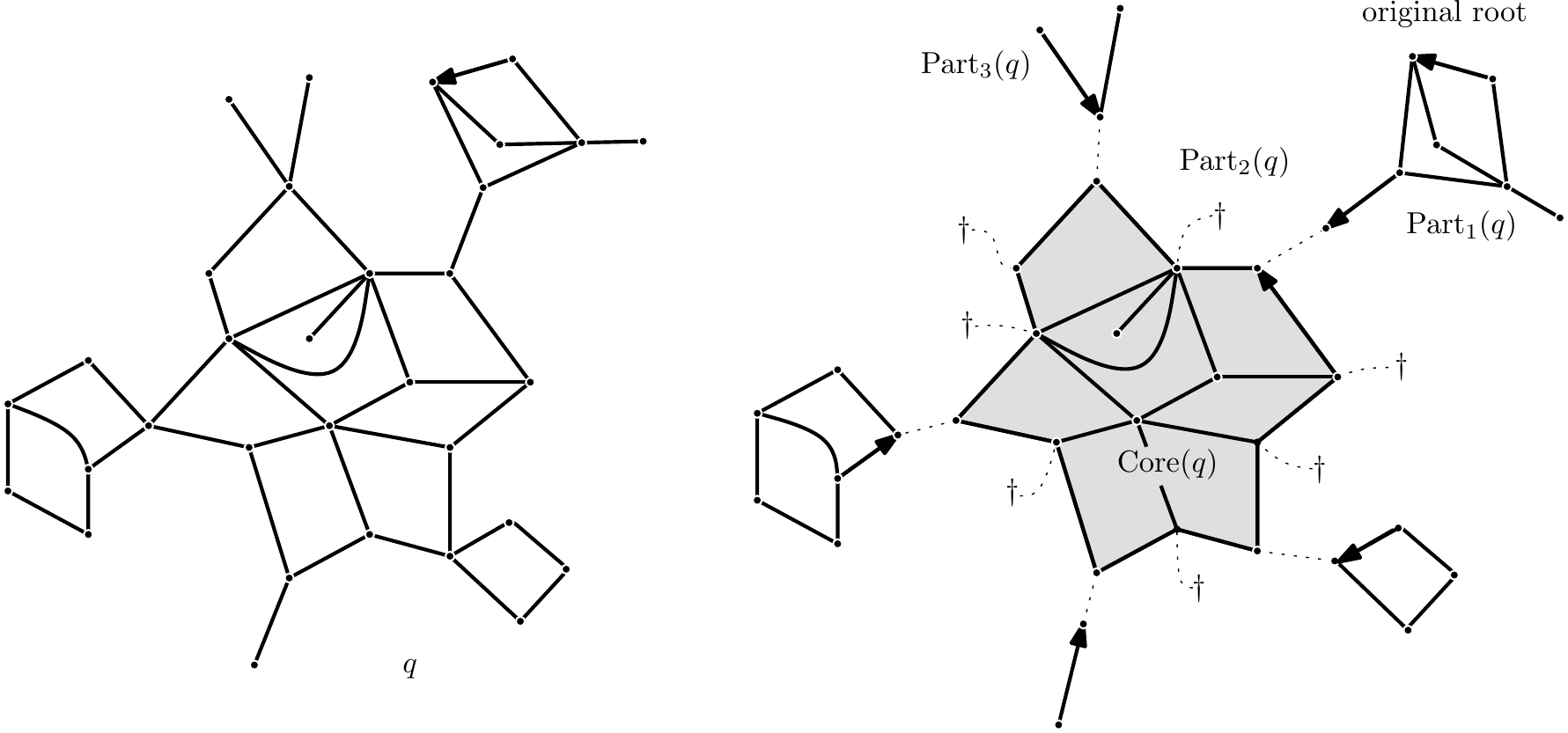} \caption{ \label{pruningfig} Illustration of the pruning.}  \end{center} \end{figure}

The quadrangulation $q$ can be recovered from $ \mathrm{Core}(q), \mathrm{Part}_1(q), \mathrm{Part}_2(q), ... , \mathrm{Part}_{2p}(q)$ if we are given a number $\mathrm{R}(q) \in \{1, ... , |\partial \mathrm{Part}_{1}(q)|+1\}$ to specify the location of the original root edge of the initial quadrangulation: On the first, second, ..., $| \partial \mathrm{Part}_{1}(q)|$-th oriented edge of $\mathrm{Part}_{1}(q)$ if $\mathrm{R}(q) \in \{1, ... , |\partial \mathrm{Part}_{1}(q)|\}$, or on the core just before $\mathrm{Part}_{1}(q)$ if $ \mathrm{R}(q) = | \partial \mathrm{Part}_{1}(q)|+1$. 

If there is no largest irreducible component, we set $ \mathrm{Core}(q)$ and all the components to be equal to $\dag$ and $ \mathrm{R}(q)= 0$ by convention. The pruning is still possible when we deal with a rooted quadrangulation with boundary that contains a unique infinite irreducible component,  which is automatically the core.

\begin{proposition}\label{oneend} 
For every $p \in \{1,2,3, ...\}$, almost surely $Q_{\infty,p}$ has 
only one infinite irreducible component.
\end{proposition}

The proof is easy using the construction of $Q_{\infty,p}$ from $ \mathbf{B}_{\infty,p}$ and the fact that $ \mathbf{B}_{\infty,p}$ contains only one infinite labeled tree.
Details are left to the reader.



Recall that $Q_{n,p}$ (resp.\,$\widetilde{{Q}}_{n,p}$) is uniformly distributed over $\mathcal{Q}_{n,p}$ (resp.\,$\widetilde{\mathcal{Q}}_{n,p}$). We also denote by $U_k$ a uniform variable over $ \{1,2,..., k\}$. Fix $p\geq 1$ and $q\geq 1$, and let $g,f_1, ... , f_{2q}$ be positive bounded continuous functions for the distance $\mathrm{d_{map}}$ and $e: \mathbb{R} \to \mathbb{R}_+$ be a bounded positive continuous function.
As an immediate consequence of Proposition \ref{oneend} and Theorem \ref{CMM10bound} we deduce that $Q_{n,p}$ has a largest irreducible component with a probability tending to $1$ as $n \to \infty$. Thus we have
\begin{eqnarray} \lefteqn{{E}\left[ g\big(\mathrm{Core}(Q_{n,p})\big) e(R(Q_{n,p}))\prod_{i=1}^{2q}f_i\big(\mathrm{Part}_i(Q_{n,p})\big)  \mathbf{1}_{|\partial \mathrm{Core}(Q_{n,p})|=2q}\right]} \nonumber \\ & \xrightarrow[n\to\infty]{}& {E}\left[ g\big(\mathrm{Core}(Q_{\infty,p})\big) e(R(Q_{\infty,p}))\prod_{i=1}^{2q}f_i\big(\mathrm{Part}_i(Q_{\infty,p})\big) \mathbf{1}_{|\partial \mathrm{Core}(Q_{\infty,p})|=2q}\right]. \label{cv1}\end{eqnarray}
Moreover, the pruning decomposition leads to 
\begin{eqnarray}
 \lefteqn{{E}\left[  g\big(\mathrm{Core}(Q_{n,p})\big) e(R(Q_{n,p}))\prod_{i=1}^{2q}f_i\big(\mathrm{Part}_i(Q_{n,p})\big)\mathbf{1}_{ |\partial\mathrm{Core}(Q_{n,p})|=2q }\right]} \nonumber \\ 
 &=& \frac{1}{q_{n,p}}\sum_{
   n \geq m \geq 0}
 \widetilde{q}_{m,q}{E}\big[g(\widetilde{Q}_{m,q})\big]\nonumber\\ &&\times
 \sum_{\begin{subarray}{c} p_{1}+ ... + p_{2q}= p-q\\
     m > n_{1}, ... ,n_{2q} \geq 0\\
     n_{1} + ... + n_{2q} = n-m \end{subarray}}  (2p_{1}+1){E}\left[e(U_{2p_1+1})\right] \prod_{i=1}^{2q}q_{n_{i},p_{i}} {E}\big[f_i(Q_{n_i,p_i})\big]. \label{decomposition} \end{eqnarray}
Similarly as in the proof of Proposition \ref{Kestenbis}, one can use Lemma \ref{utile} to deduce that the probability that one of the components $ \mathrm{Part}_{1}(Q_{n,p}), \mathrm{Part}_{2}(Q_{n,p}), ...$ has a size larger than $a >0$ is bounded above by $c a^{-3/2}$ for some constant $c>0$ uniformly in $n$. So we can let  $n \to \infty$ followed by $a \to \infty$ in the formula \eqref{decomposition} and obtain by \eqref{def:UIPQs},  \eqref{cv1} and asymptotics \eqref{asympqnp} and \eqref{asympqnpt} that

\begin{eqnarray}
\lefteqn{{E}\left[ g\big(\mathrm{Core}(Q_{\infty,p})\big) e(R(Q_{\infty,p}))\prod_{i=1}^{2q}f_i\big(\mathrm{Part}_i(Q_{\infty,p})\big)\mathbf{1}_{|\partial \mathrm{Core}(Q_{\infty,p})|=2q}\right] }\nonumber \\
 &=& \frac{\widetilde{C}_q }{C_p}{E}\big[g(\widetilde{Q}_{\infty,q})\big]\sum_{\begin{subarray}{c}
p_{1}+ ... + p_{2q}= p-q\end{subarray}}  (2p_{1}+1){E}\left[e(U_{2p_1+1})\right] \nonumber\\
&&\times \prod_{i=1}^{2q} \sum_{n_i=0}^{\infty}12^{-n_i}q_{n_{i},p_{i}} {E}\big[f_i(Q_{n_i,p_i})\big]. \label{decompositioninfty} \end{eqnarray} 

The last expression is the fundamental ``pruning formula''. It can be used to derive the distribution of the core and the components of $Q_{\infty,p}$ for a fixed $p$. Let us proceed.  Fix $p \geq q\geq 1$, so that \eqref{decompositioninfty} specializes when $e=f_1=\ldots=f_{2q}=1$ to
\begin{eqnarray*}
E[g(\mathrm{Core}(Q_{\infty,p}))]&=&\frac{\widetilde{C}_q}{C_{p}}E[g(\widetilde{Q}_{\infty,q})]
\sum_{p_1+\ldots+p_{2q}=p-q}(2p_1+1)\prod_{i=1}^{2q}\sum_{n_i=0}^\infty12^{-n_i}q_{n_i,p_i}\\
&=&\frac{p}{q} \frac{\widetilde{C}_q}{C_{p}}E[g(\widetilde{Q}_{\infty,q})]
\sum_{p_1+\ldots+p_{2q}=p-q}\prod_{i=1}^{2q}\sum_{n_i=0}^\infty12^{-n_i}q_{n_i,p_i}\, ,
\end{eqnarray*}
where we got rid of the term $(2p_1+1)$ by an obvious symmetry argument, to the cost of adding the prefactor $(2q)^{-1}\sum_{i=1}^{2q}(2p_i+1)=p/q$. Recalling the definition of the bivariate function $W(g,z)$ and $W_{c}(z)=W(12^{-1},z)$, we can further re-write the last expression as 
\begin{eqnarray*}
\lefteqn{\frac{p}{q}\frac{\widetilde{C}_q}{C_{p}} E[g(\widetilde{Q}_{\infty,q})] 8^{p-q}W_{c}(8^{-1})^{2q}
\sum_{p_1+\ldots+p_{2q}=p-q}\prod_{i=1}^{2q}\frac{\sum_{n_i=0}^\infty12^{-n_i}8^{-p_i}q_{n_i,p_i}}{W_{c}(8^{-1})}}\\
&=&\frac{p}{q}\frac{(9/2)^{-q}\widetilde{C}_q}{8^{-p}C_{p}} E[g(\widetilde{Q}_{\infty,q})]
\sum_{p_1+\ldots+p_{2q}=p-q}\prod_{i=1}^{2q}\frac{\sum_{n_i=0}^\infty12^{-n_i}8^{-p_i}q_{n_i,p_i}}{W_{c}(8^{-1})}\, ,
\end{eqnarray*}
where we used the fact that $W_{c}(8^{-1})=4/3$, see \eqref{WWW}.
We interpret the last sum as 
$P(Z_1+Z_2+\ldots+Z_{2q}=p-q)$, where $Z_1,\ldots,Z_{2q}$ are independent random variables with common distribution 
$$P(Z_1=r)=\frac{\sum_{n=0}^\infty12^{-n}8^{-r}q_{n,r}}{W_{c}(8^{-1})}\, ,\qquad r\geq 0\, .$$
We just proved \begin{theorem}[Pruning with fixed perimeter] \label{thm:fixed}For every $p\geq q \geq 1$, conditionally on  the event $\{|\partial \mathrm{Core}(Q_{\infty,p})|=2q\}$ of probability
 \begin{eqnarray*}P( |\partial \mathrm{Core}(Q_{\infty,p})|=2q) &=& \frac{q^{-1}(9/2)^{-q} \widetilde{C}_{q}}{p^{-1}8^{-p} C_{p} } P( Z_{1} + ... + Z_{2q}= p-q),  \end{eqnarray*}  the core of $Q_{\infty,p}$ is distributed as a simple boundary UIPQ with perimeter $2q$.\end{theorem}

\subsection{Proof of Theorem \ref{main}}
\label{sec:proof-theor-refm}

As a first application of Theorem \ref{thm:fixed}, let us now prove Theorem \ref{main}.
To this end, we first make some preliminary observations. 
By definition of $W(g,z)$, the generating function of $Z_1$ is given by 
\begin{eqnarray*}
E[s^{Z_1}]&=&\frac{W_c(s8^{-1})}{W_c(8^{-1})}=\frac{(1-\sqrt{1-s})(2s-1+\sqrt{1-s})}{s^2}\\
&\underset{s\uparrow 1}{=}&s+2(1-s)^{3/2}+o((1-s)^{3/2})\, .
\end{eqnarray*}
By standard results on stable domains of attraction \cite{BGT89}, this expression entails that the random variable $Z_1$ is in the domain of attraction of a stable random variable with exponent $3/2$. More precisely, since moreover $E[Z_1]=1$ by differentiating the previous expression, it holds that
$$\frac{Z_1+\ldots+Z_n-n}{n^{2/3}}\xrightarrow[n\to\infty]{(d)}\mathcal{Z}'\, ,$$
where the Laplace transform of $\mathcal{Z}'$ is given by $E[\exp(-\lambda \mathcal{Z}')]=\exp(2\lambda^{3/2})$ for every $\lambda\geq 0$.  More precisely, if $h$ denotes the density of the law of $\mathcal{Z}'$, then the Gnedenko-Kolmogorov local limit theorem for lattice variables (see \cite[Theorem 4.2.1]{IL71}) entails that
\begin{equation}\label{eq:1}
\sup_{k\in\Z}\left|n^{2/3}P(Z_1+Z_2+\ldots+Z_{n}=n+k)-h\Big(\frac{k}{n^{2/3}}\Big)\right|\underset{n\to\infty}{\longrightarrow}0\, .
\end{equation}

Now, for a given $q\geq 1$, set $p=3q$ in equation \eqref{decompositioninfty}. By using \eqref{eq:1} with $n=2q$ and $k=0$, we obtain $P(Z_1+\ldots+Z_{2q}=2q) \sim h(0) (2q)^{-2/3}$ as $q \to \infty$.  On the other hand, the asymptotic behavior for $C_p,\widetilde{C}_q$ entails that (still when $p=3q$)
 \begin{eqnarray}\label{eq:2}
\frac{q^{-1}(9/2)^{-q}\widetilde{C}_q}{p^{-1}8^{-p}C_p}&\xrightarrow[q\to\infty]{}&3\, .
 \end{eqnarray}
From this and Theorem \ref{thm:fixed}, we conclude that $$P(|\partial \mathrm{Core}(Q_{\infty,3q})| =2q) \quad \sim \quad A q^{-2/3}\, ,$$ where $A=3\cdot2^{-2/3}\cdot h(0)\in (0,\infty)$. By the first assertion of Theorem \ref{thm:fixed}, we deduce that there exists $A'\in (0,\infty)$ such that for any non-negative measurable function $g$,
$$E[g(\widetilde{Q}_{\infty,q})] = E\big[g(\mathrm{Core}(Q_{\infty,3q})) \,\big| \,|\partial	\mathrm{Core}(Q_{\infty,3q}) | =2q\big]\leq A' q^{2/3}E[g(\mathrm{Core}(Q_{\infty,3q}))]\, .$$
From this, Theorem \ref{diamext}, and the obvious fact that $\mathrm{aper}(\mathrm{Core}(Q_{\infty,3q}))\leq \mathrm{aper}(Q_{\infty,3q})$, we conclude that \begin{eqnarray*}
  P(\mathrm{aper}(\widetilde{Q}_{\infty,q})\geq \lambda\sqrt{q}) &\leq &
A' q^{2/3} P(\mathrm{aper}(\mathrm{Core}(Q_{\infty,3q}))\geq \lambda\sqrt{q})\\
&\leq & A' q^{2/3}  P(\mathrm{aper}(Q_{\infty,3q})\geq \lambda\sqrt{q})\\
&\leq & A'c_1 q^{2/3} \exp(-3^{-1/3}c_2\lambda^{2/3})\, .
\end{eqnarray*}
This yields Theorem \ref{main}. \endproof 

\begin{remark} Theorem \ref{diamext} entails that $(p^{-1/2} \mathrm{aper}(Q_{\infty,p}))_{p \geq 0}$ is tight. We believe that a similar property holds for $({p}^{-1/2} \mathrm{aper}(\widetilde{Q}_{\infty,p}))_{p \geq 0}$, but this is not a direct consequence of our results.  In fact, we believe that  $({p}^{-1/2} \mathrm{aper}(Q_{\infty,p}))_{p \geq 0}$ and $({(3p)}^{-1/2} \mathrm{aper}(\widetilde{Q}_{\infty,p}))_{p \geq 0}$ converge in distribution to the same non-degenerate random variable.  \end{remark}

\subsection{Asymptotics of the perimeters}

As a second application of Theorem \ref{thm:fixed}, we will see that as $p \to \infty$ the core of $Q_{\infty, p}$ has a perimeter which is roughly a third of the original quadrangulation. This supports \cite[Section 5]{BG09} where the authors proved that quadrangulations with simple boundary of perimeter $p$ have the same large scale structure as quadrangulations with general boundary of perimeter $3p$. 

\begin{proposition} \label{1/3} We have  the following convergence in probability\begin{eqnarray*} \frac{|\partial \mathrm{Core}( Q_{\infty,p})|}{2p} & \xrightarrow[p\to\infty]{(P)}& \frac{1}{3}.  \end{eqnarray*}
More precisely, it holds that 
 \begin{eqnarray*}\frac{|\partial \mathrm{Core}( Q_{\infty,p})|-2p/3}{p^{2/3}} & \xrightarrow[p\to\infty]{(d)} &\mathcal{Z}\, ,  \end{eqnarray*}
where $\mathcal{Z}$ is a spectrally negative stable random variable with exponent $3/2$, with Laplace transform given by 
$$E[\exp(\lambda \mathcal{Z})]=\exp\bigg(\Big(\frac{2}{3}\Big)^{5/2}\lambda^{3/2}\bigg)\, ,\qquad \lambda\geq 0\, .$$
\end{proposition}

\proof 
Let $p\geq 1$ and $x\in \R$, and let $q=\lfloor p/3\rfloor+\lfloor xp^{2/3}\rfloor$. We again use the local limit theorem \eqref{eq:1} by specializing it to $n=2q$ and $k=p-3q$, and utilize the asymptotic equivalents for $C_{p},\widetilde{C}_q$ with the same limit as in \eqref{eq:2}. Together with Theorem \ref{thm:fixed} this implies
$$p^{2/3}P\Big({\frac{1}{2}|\partial\mathrm{Core}(Q_{\infty,p})|-\lfloor p/3\rfloor}=\lfloor xp^{2/3}\rfloor\Big)\xrightarrow[p\to\infty]{}\frac{3^{5/3}}{2^{2/3}}h\Big(-\frac{3^{5/3}}{2^{2/3}}x\Big)\, .$$
By Scheff\'e's lemma and elementary computations using the Laplace transform of $h$, this implies the claim on convergence in distribution in the statement.  The first claim on convergence in probability is a simple consequence of the latter.
\endproof

\subsection{Randomizing the perimeters}\label{sec:rand-bound-length}

In this section, we argue that \eqref{decompositioninfty} gives a particularly nice probabilistic interpretation of the pruning operation, to the cost of randomizing the perimeters of the maps under consideration. 
Let us introduce some notation.  Let $C(z)=\sum C_p z^p$ be the generating function of the $C_p$'s (with $C_0=0$) and set $\vec{W}_c(z) := W_{c}(z) + 2z \partial_{z}W_{c}(z)$. Notice that $\sum_{n\geq 0} q_{n,p}(2p+1)12^{-n}= [z^p](W_{c}(z) + 2z \partial_{z}W_{c}(z))=[z^p]\vec{W}_c(z)$. From the exact expressions of $C_{p}(.)$ \eqref{asympqnp} 
 and $W(.,.)$ \eqref{WWW} we get that
\begin{eqnarray} C(z) &=& \frac{2 z}{\sqrt{\pi} (1 - 8 z)^{3/2}}, \nonumber \\
W_c(z) &=& \frac{(1-8z)^{3/2}-1+12z}{24z^2}, \nonumber \\
\vec{W}_c(z) &=& \frac{1 - 4 z-\sqrt{1-8z}}{8 z^2} \label{look}.
\end{eqnarray}
For $0 \leq  z \leq 1/8$, we denote by $Q_{f,z}$ (resp.\,$\vec{Q}_{f,z}$) a random finite quadrangulation with general boundary such that the size of $Q_{f,z}$ equals $n$ and its perimeter $2p$ with probability $12^{-n}z^{p}W_c(z)^{-1}$ (resp.\,$(2p+1)12^{-n}z^{p}\vec{W}_c(z)^{-1}$). Since $W_c(1/8)$ and $\vec{W}_c(1/8)$ are finite, both $Q_{f,z}$ and $\vec{Q}_{f,z}$ make sense for $z=1/8$. We call these random quadrangulations ``free critical Boltzmann (extra rooted) quadrangulations with parameter $z$". Finally, for $0 < z < 1/8$, let $ \mathcal{P}_z$ and $ \widetilde{ \mathcal{P}}_z$ be  random variables distributed according to 

 \begin{eqnarray}& \displaystyle {P}( \mathcal{P}_z=p) \  =\ \frac{z^p C_p}{C(z)}, \quad\quad  \mbox{and} \quad \quad  {P}( \widetilde{ \mathcal{P}}_z=q) \  =\   \widetilde{C}_q \big(zW^2_c(z)\big)^q \frac{\vec{W}_{c}(z)}{C(z)W_{c}(z)}.& \label{defpwp}  \end{eqnarray}  Where we recall that $C_0= \widetilde{C}_0=0$. The lines leading to \eqref{decompositioninfty} (or a direct calculation) show that $ \widetilde{C}(z W_{c}^2(z))\vec{W}_{c}(z) = C(z)W_{c}(z)$ for every $0 \leq z < 1/8$ (where $ \widetilde{C}$ is the generating function of the $ \widetilde{C}_{q}$'s) so that $ \widetilde{\mathcal{P}}_{z}$ is well-defined. 
 
 In the remaining of this work if $\mathcal{P}$ is a integer-valued random variable, we denote by $Q_{\infty,\mathcal{P}}$  a random variable such that conditionally on $\{ \mathcal{P}=p\}$, $Q_{\infty, \mathcal{P}}$ is distributed as $Q_{\infty,p}$ (and similarly of the ``$\sim$'' analog).
 \medskip 
 
Fix $ z \in (0 , 1/8)$. Now that the reader is acquainted with this notation we multiply both members of \eqref{decompositioninfty} by $z^pC_p$, sum over all $p \geq 1$ and divide after-all by $C(z)$ to  deduce that 

\begin{eqnarray}
\lefteqn{ {E}\left[ g\big(\mathrm{Core}(Q_{\infty,\mathcal{P}_z})\big) e(R(Q_{\infty,{ \mathcal{P}}_z}))\prod_{i=1}^{2q}f_i\big(\mathrm{Part}_i(Q_{\infty,\mathcal{P}_z})\big)\mathbf{1}_{|\partial \mathrm{Core}(Q_{\infty,{ \mathcal{P}}_z})|=2q}\right] }\nonumber \\
& =&\frac{\vec{W}_{c}(z)}{C(z)W_{c}(z)} \widetilde{C}_q \big(zW^2_c(z)\big)^q \left ({E}\big[g(\widetilde{Q}_{\infty,q})\big]  {E}\Big[f_1(\vec{Q}_{f,z}) {E}\left[e\big(U_{|\partial \vec{Q}_{f,z}|+1}\big)\right]\Big]\prod_{i=2}^{2q} {E}\big[f_i(Q_{f,z})\big]\right)  \nonumber \\
 &=& P(\widetilde{\mathcal{P}}_z=q)\left( E\left[ g(\widetilde{Q}_{\infty,q} )\right] {E}\Big[f_1(\vec{Q}_{f,z}) {E}[e(U_{|\partial \vec{Q}_{f,z}|+1})]\Big]\prod_{i=2}^{2q} {E}\big[f_i(Q_{f,z})\big]\right). \nonumber\\ &&\label{decompositioninftyboltz} \end{eqnarray}
Thus we proved: 

\begin{theorem}[Pruning with random perimeter]\label{thmpruning} For every $ z\in (0, 1/8)$, we have the following equality in distribution \begin{eqnarray*}\mathrm{Core}\big(Q_{\infty, \mathcal{P}_z}\big) & \overset{(d)}{=}& \widetilde{Q}_{\infty,\widetilde{\mathcal{P}}_z}. \end{eqnarray*} Furthermore, conditionally on $|\partial \mathrm{Core}(Q_{\infty, \mathcal{P}_z})|$, the core and the components of $Q_{\infty, \mathcal{P}_z}$ are independent the latter being distributed as follows: the first component is distributed according to $\vec{Q}_{f,z}$ (the location $R(Q_{\infty, \mathcal{P}_z})$ of the root being uniform over $\{1,2, ... , | \partial \mathrm{Part}_1(Q_{\infty,p})|+1\}$), and all the other $|\partial \mathrm{Core}(Q_{\infty, \mathcal{P}_z})|-1$ components are distributed according to $Q_{f,z}$.  \end{theorem}

Using the exact expression of $C(.)$ we deduce that for $z \in (0,1/8)$ we have
 \begin{eqnarray*}E[ \mathcal{P}_z] &=& \frac{z C'(z)}{C(z)} \quad =\quad  \frac{1+4z}{1-8z}.  \end{eqnarray*}
Thus the average perimeter of $ Q_{\infty, \mathcal{P}_z}$ is asymptotically equivalent to $ 3(1-8z)^{-1}$ as $z \to 1/8$ and it is easy to see from the singularity analysis of $C(.)$ that $ \mathcal{P}_z \to \infty$ in distribution as $z \to 1/8$. Using Proposition \ref{1/3} (or by a direct analysis) we see that $ \widetilde{ \mathcal{P}}_z \to \infty$ as well as $z \to 1/8$.


\subsection{Interpretation of labels when $Q_{\infty,p}= \Phi( \mathbf{B}_{\infty,p})$}
In this section we finally prove how the pruning can be used to deduce the second part of Theorem \ref{CMM10bound} from the results of \cite{CMMinfini}. Let us recall the setting. Fix $ p \in \{1,2,3,...\}$ and let $ \mathbf{B}$ be a uniform infinite labeled treed bridge of length $p$. We consider $Q_{\infty,p} = \Phi( \mathbf{B})$.  We aim at showing that, with the identification of the vertices of $Q_{\infty,p}$ with those of $\mathbf{B}$, a.s.\, for every $u,v \in Q_{\infty,p}$ we have
 \begin{eqnarray*}\lim_{z \to \infty} \big(\mathrm{d_{gr}}(u,z)-\mathrm{d_{gr}}(v,z)\big) &=& \ell_{\mathbf{B}}(u)-\ell_{\mathbf{B}}(v).   \end{eqnarray*}
\proof[Proof of the second part of Theorem \ref{CMM10bound}] 
In order to prove the last display, it only suffices to prove that the left-hand side actually a.s.\,has a limit as $z \to \infty$: Let us call this fact the property $(*)$. Then the last display follows from an adaptation of the end of the proof of \cite[Lemma 5]{CMMinfini}: For $u,v \in \mathbf{B}$ we can consider two simple geodesics $\gamma_{u}$ and $\gamma_{v}$ starting from any corners associated with $u$ and $v$ in $\mathbf{B}$. These two paths eventually merge. If we assume that property $(*)$ holds then we take $z \to \infty$ along the geodesics after the merging point: For such $z$ we have  $\mathrm{d_{gr}}(u,z)-\mathrm{d_{gr}}(v,z) = \ell_ \mathbf{B}(u)-\ell_ \mathbf{B}(v)$, which proves the claim.\\
It thus suffices to prove that $(*)$ holds for $Q_{\infty,p}$. In fact, it is easy to see that we can restrict our attention to those $u,v$ that belong to the Core of $Q_{\infty,p}$. Thanks to Theorem \ref{thmpruning}, for any $1 \leq q \leq p$, conditionally on $\{\partial \mathrm{Core}(Q_{\infty,p})=2q\}$ we have $ \mathrm{Core}(Q_{\infty,p}) = \widetilde{Q}_{\infty,q}$ in distribution. We are thus reduced to prove that $\widetilde{Q}_{\infty,q}$ satisfies $(*)$ for all $q \geq 1$. To show this, we will make some plastic surgery with $\widetilde{Q}_{\infty,p}$ in order to come back to the setup of \cite{CMMinfini} which deals with the full-plane UIPQ.

\begin{figure}[!h]
 \begin{center}
 \includegraphics[width=14cm]{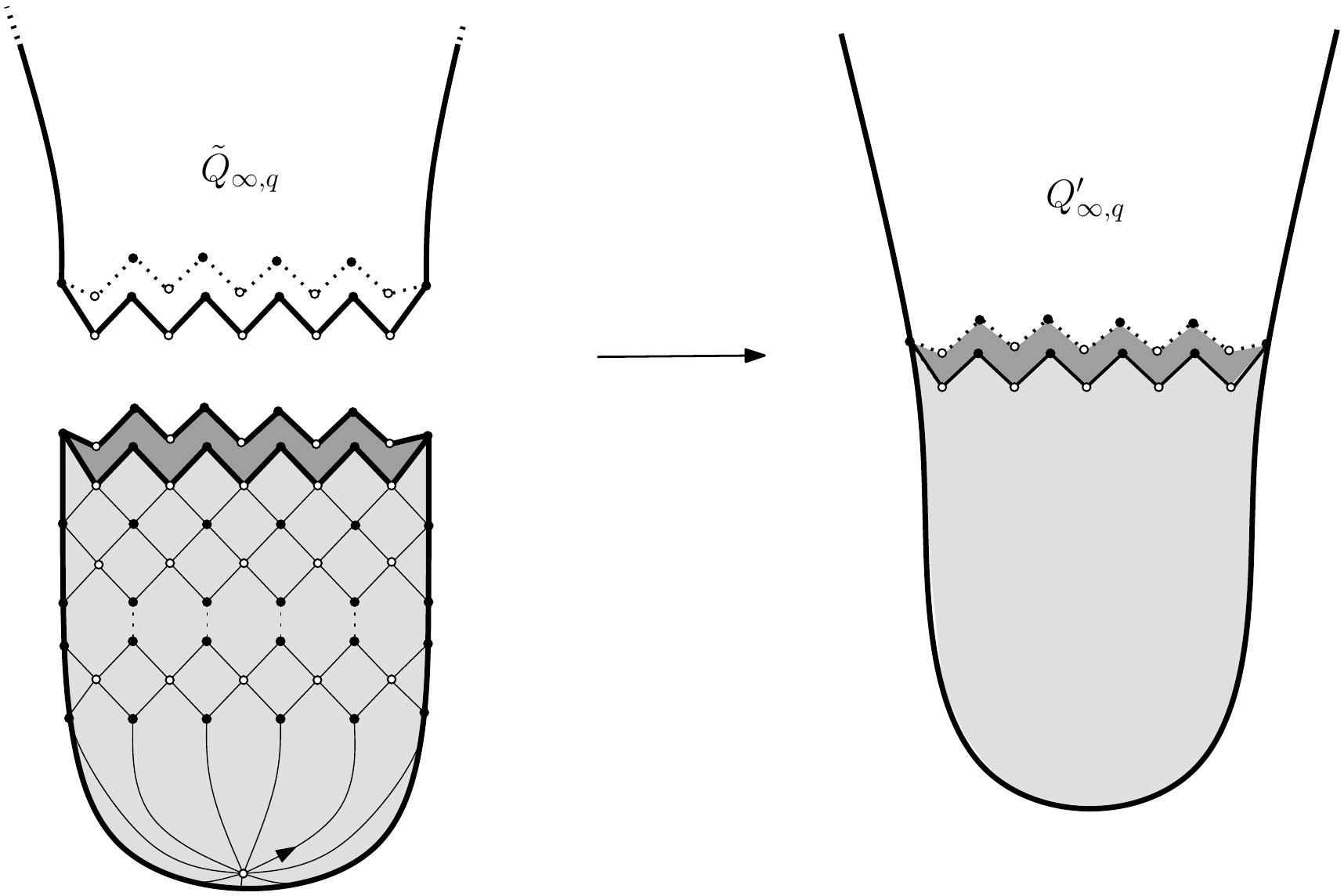}
 \caption{ \label{chirurgie}The filling operation.}
 \end{center}
 \end{figure}

More precisely, let us consider the quadrangulation ${Q}_{\infty,q}'$ obtained after filling the external face of $\widetilde{Q}_{\infty,q}$ with a quadrangulation with a simple boundary of perimeter $2q$ made of $2q+1$ ``layers'' of quadrangles such that the last layer is connected to the root edge, see Fig.\,\ref{chirurgie}.

 This operation is reversible, that is, given ${Q}'_{\infty,q}$ and the number $q$, we can recover $\widetilde{Q}_{\infty,p}$. Most importantly, the filing has been done in such a way that
 for any $u,v \in \widetilde{Q}_{\infty,q}$ we have 
 \begin{eqnarray} \label{clef} \mathrm{d}_{ \mathrm{gr}}^{{Q}_{\infty,q}'}(u,v) & =& \mathrm{d}_{ \mathrm{gr}}^{\widetilde{Q}_{\infty,q}}(u,v). \end{eqnarray} Indeed it is easy to see that for any $u,v \in \widetilde{Q}_{\infty,q}$ we can find a geodesic path between $u$ and $v$ that does not enter the grated region. Furthermore, by the spatial Markov property of the UIPQ (see \cite{BCsubdiffusive}) we deduce that the law of ${Q}_{\infty,q}'$ is absolutely continuous with respect to the law of the UIPQ.  Since the UIPQ almost surely satisfies the property $(*)$ (see \cite{CMMinfini}) we deduce that ${Q}_{\infty,q}'$ and by \eqref{clef} that $\widetilde{Q}_{\infty,q}$ also satisfies it a.s. \endproof
 
 \begin{remark} This surgical operation can also be used to transfer other ``ergodic'' properties of the standard UIPQ towards UIPQ with boundaries. \end{remark}

\section{Open boundary, open questions}

\subsection{UIPQ with infinite boundary}
In this section we let $p \to \infty$ and define the UIPQ with infinite general and simple boundary of infinite perimeter. 
We then extend the pruning procedure to these infinite quadrangulations. The proofs are only sketched or left to the reader.

\subsubsection{General boundary} We start by introducing the limit of the uniform treed bridges as $p \to \infty$.\\

Let $(X_n)_{n\in \mathbb{Z}}$ be a two-sided simple random walk starting from $0$ at $0$ and having uniform increments in $\{+1,-1\}$. Independently of $(X_n)_{n\in\mathbb{Z}}$, let $(\Theta_i)_{i \in \mathbb{Z}}$ be a sequence of independent uniform labeled geometric critical Galton-Watson trees. The object $\mathbf{B}_{\infty,\infty}=((X_n)_{n\in \mathbb{Z}}, (\Theta_i)_{i \in \mathbb{Z}})$ is \emph{the uniform infinite treed bridge of infinite length}. It obviously appears as a limit of the uniform infinite treed bridge of length $2p$ as $p \to \infty$ in the following sense: Let $\mathbf{B}_{\infty,p}=(B_p; \Theta_{1,p}, ... , \Theta_{p,p})$ be a uniform infinite treed bridge of length $2p$ with  $B_p=(X_{1}^{(p)}, ... , X^{(p)}_{2p})$. Then for any $m \geq 1$ we have  \begin{eqnarray*} (X_{[i]}^{(p)})_{ -m \leq i \leq m}& \xrightarrow[p\to\infty]{(d)} & (X_i)_{-m\leq i \leq m}, \end{eqnarray*}  
when $[i]$ stands for the representative of $i$ modulo $2p$ that belongs to $\{1,...,2p\}$. Furthermore the trees grated on the down-steps $[i]$ such that $-m \leq i \leq m$ asymptotically are i.i.d. critical geometric GW trees since the probability that one of these trees is the infinite one tends to $0$ as $p \to \infty$.\\ 

\begin{figure}[!h]
 \begin{center}
 \includegraphics[width=16cm]{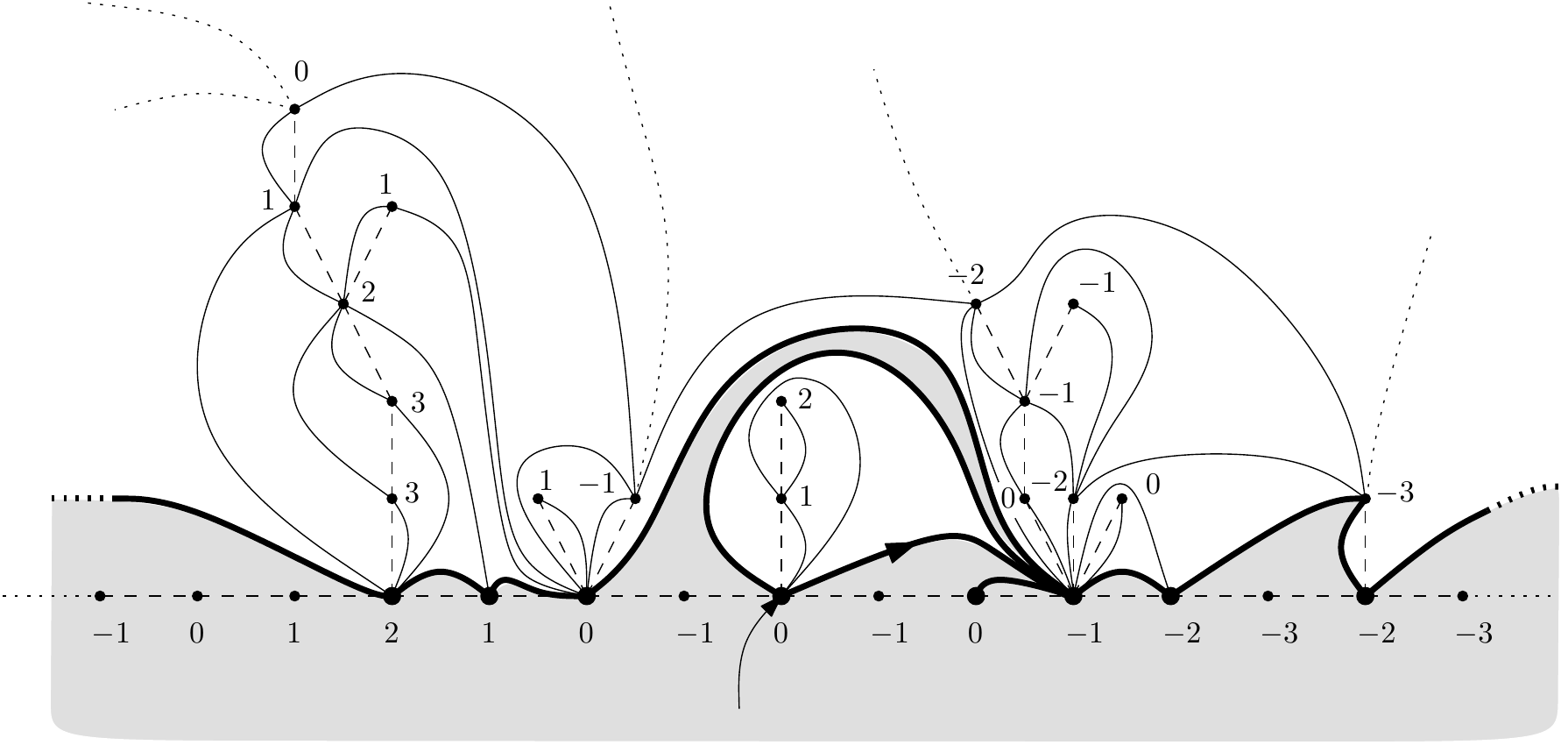}
 \caption{ \label{htuipq}Second extension of $\Phi$.}
 \end{center}
 \end{figure}

 We can associate with the infinite treed bridge $ \mathbf{B}_{\infty,\infty}$ a representation by grafting the trees $\Theta_i$ to the down-steps of the walk $(X_n)$, see Fig.\,\ref{htuipq}. We then (once again) extend the Schaeffer mapping $\Phi$ to this object in a straightforward manner and define   a random infinite quadrangulation with an infinite perimeter denoted by $\Phi( \mathbf{B}_{\infty,\infty})$, see Fig.\,\ref{htuipq}. 

\begin{proposition} \label{half} We have the following convergence 
\begin{eqnarray*}
Q_{\infty,p} & \xrightarrow[p\to\infty]{(d)}& Q_{\infty,\infty}, 
\end{eqnarray*}
in distribution for $ \mathrm{d_{map}}$, where $Q_{\infty,\infty}$ 
is random infinite rooted quadrangulation with an infinite boundary which can be constructed from the uniform infinite treed bridge of infinite perimeter via the extended Schaeffer mapping, that is $ Q_{\infty,\infty} = \Phi( \mathbf{B}_{\infty,\infty})$ in distribution.
\end{proposition}
This is essentially an adaptation of the first part of Theorem \ref{CMM10bound} (and Proposition \ref{continuity}). The proof is left to the interested reader.   \bigskip

\subsubsection{Simple boundary.} Recall that $ \mathcal{P}_z \to \infty$ in distribution as $ z \uparrow 1/8$. Thus it follows from the last theorem that $ \mathcal{Q}_{\infty, \mathcal{P}_z} \to Q_{\infty,\infty}$ as $z \to 1/8$. Using Theorem \ref{thmpruning} we deduce that $ \widetilde{Q}_{\infty, \widetilde{\mathcal{P}}_z}$ converge towards some random infinite quadrangulation with an infinite simple boundary as $z \to 1/8$. We denote this limit by $\widetilde{Q}_{\infty,\infty}$ and call it the UIPQ with infinite simple boundary or UIPQ of the half-plane.  

Angel \cite{Ang05} defined and studied the analog of $\widetilde{Q}_{\infty,\infty}$ in the triangulation case. His approach can be adapted to the quadrangulation case to show that 
 \begin{eqnarray*}\widetilde{Q}_{\infty,p} & \xrightarrow[p\to\infty]{(d)} & \widetilde{Q}_{\infty,\infty}, \end{eqnarray*} for $ \mathrm{d_{map}}$. One of the advantages of working with such objects is the very simple form that takes the spatial Markov property, see \cite{Ang05,ACpercopeel}.

The pruning procedure can also be extended to $Q_{\infty,\infty}$ (one can show that $Q_{\infty,\infty}$ has only one infinite irreducible component almost surely). The following statement can be seen as an extension to $z=1/8$ of Theorem \ref{thmpruning}. 
\begin{proposition}
The core and the components of $Q_{\infty,\infty}$ are all independent, the core being distributed as $\widetilde{Q}_{\infty,\infty}$, the first component as $\vec{Q}_{f,1/8}$, and the other components as $Q_{f, 1/8}$.
\end{proposition}

\subsection{Comments, questions}
Extending the techniques of this paper, it is possible to study a variant of the aperture in the case of $Q_{\infty,\infty}$ and translate the results to the simple boundary case via the pruning procedure extended in Theorem \ref{half}. We present here a couple of open questions related to the models $Q_{\infty,\infty}$ and $ \widetilde{Q}_{\infty,\infty}$.

\begin{open} In the construction $Q_{\infty,\infty} = \Phi( \mathbf{B}_{\infty,\infty})$, is it the case that the $\ell_{ \mathbf{B}_{\infty,\infty}}$-labels have the same interpretation as in Theorem \ref{CMM10bound}, that is, for every $u,v \in Q_{\infty,\infty}$ we have 
 \begin{eqnarray*}  \lim_{z \to \infty} \big( \mathrm{d_{gr}}(z,u)- \mathrm{d_{gr}}(z,v)\big) &=& \ell_{ \mathbf{B}_{\infty,\infty}}(u)-\ell_{ \mathbf{B}_{\infty,\infty}}(v).\end{eqnarray*}
 \end{open}
 
 We next move to surgical considerations. Consider two copies of $ \widetilde{Q}_{\infty,\infty}$ and glue them together along the boundary with coinciding roots to form a rooted quadrangulation of the plane denoted by $ \mathfrak{Q}_{\infty}$.  We claim that the law of $ \mathfrak{Q}_{\infty}$ is singular with respect to the law of $Q_{\infty}$. Indeed, it is easy to construct two infinite (simple) geodesics in $ \mathfrak{Q}_{\infty}$ starting from the origin and that are eventually non intersecting. However, two such geodesics do not exist in the case of the UIPQ, see \cite{CMMinfini}. Consequently, the full-plane UIPQ is \emph{not} the result of the gluing of two independent half-plane UIPQ with simple boundary. A more interesting gluing is the following:

 \begin{open} Consider the ``closing'' operation that consists in zipping the boundary of $\widetilde{Q}_{\infty,\infty}$ to get an infinite rooted quadrangulation $ \mathscr{Q}_{\infty}$ with an infinite self-avoiding path on it. Is the law of $ \mathscr{Q}_{\infty}$ absolutely continuous or singular with respect to the law of $Q_{\infty}$? Study the self-avoiding walk obtained on it: In particular, is it diffusive? \end{open} 


\bibliographystyle{alpha}

\end{document}